\documentclass[11pt]{article}
\usepackage{amsfonts}
\usepackage[leqno]{amsmath}
\usepackage{graphicx}
\usepackage{hyperref}
\hypersetup{
	colorlinks=true,
	linkcolor=blue,
	filecolor=magenta,      
	urlcolor=blue,
	citecolor=red
}
\usepackage{latexsym}
\usepackage{amsmath,amsfonts,amssymb,amsthm,mathrsfs,euscript,makeidx,color}
\usepackage{enumerate}
\allowdisplaybreaks[1]
\usepackage{multirow}

\def\sqr#1#2{{\vcenter{\vbox{\hrule height.#2pt
				\hbox{\vrule width.#2pt height#1pt \kern#1pt \vrule width.#2pt}
				\hrule height.#2pt}}}}
\def\signed #1{{\unskip\nobreak\hfil\penalty50
		\hskip2em\hbox{}\nobreak\hfil#1
		\parfillskip=0pt \finalhyphendemerits=0 \par}}
\def\endpf{\signed {$\sqr69$}}

\def\3n{\negthinspace \negthinspace \negthinspace }
\def\2n{\negthinspace \negthinspace }
\def\1n{\negthinspace }
\def\bel{\begin{equation}\label}
	\def\eel{\end{equation}}

\def\dbE{\mathbb{E}}
\def\dbF{\mathbb{F}}

\def\dbH{\mathbb{H}}

\def\dbP{\mathbb{P}}

\def\dbR{\mathbb{R}}
\def\dbS{\mathbb{S}}
\def\dbT{\mathbb{T}}

\def\sU{\mathscr{U}}

\def\={\buildrel \triangle \over =}

\def\ds{\displaystyle}

\def\ns{\noalign{\ss}}
\def\d{\delta}
\def\e{\varepsilon}

\def\l{\lambda}
\def\m{\mu}
\def\n{\nu}
\def\si{\sigma}

\def\f{\varphi}
\def\th{\theta}

\def\i{\infty}
\def\D{\Delta}
\def\Th{\Theta}
\def\L{\Lambda}

\def\F{\Phi}
\def\O{\Omega}

\def\cB{{\cal B}}

\def\cF{{\cal F}}

\def\cJ{{\cal J}}

\def\cL{{\cal L}}

\def\cP{{\cal P}}

\def\cS{{\cal S}}

\def\cW{{\cal W}}

\def\BP{{\bf P}}

\def\BX{{\bf X}}
\def\BY{{\bf Y}}
\def\BZ{{\bf Z}}

\def\Bp{{\bf p}}

\def\Bu{{\bf u}}

\def\D{\Delta}
\def\Th{\Theta}
\def\L{\Lambda}

\def\F{\Phi}
\def\O{\Omega}

\def\ss{\smallskip}
\def\ms{\medskip}
\def\bs{\bigskip}
\def\q{\quad}
\def\qq{\qquad}
\def\hb{\hbox}

\def\limsup{\mathop{\overline{\rm lim}}}

\def\lan{{\langle}}
\def\ran{{\rangle}}

\def\h{\widehat}
\def\wt{\widetilde}

\def\cd{\cdot}

\def\les{\leqslant}
\def\ges{\geqslant}

\def\({\Big (}
\def\){\Big )}
\def\[{\Big[}
\def\]{\Big]}
\def\lan{\langle}
\def\ran{\rangle}

\def\bde{\begin{definition}\label}
	\def\ede{\end{definition}}

	\def\bel{\begin{equation}\label}
		\def\ee{\end{equation}}
	\def\bt{\begin{theorem}\label}
		\def\et{\end{theorem}}
	\def\bc{\begin{corollary}\label}
		\def\ec{\end{corollary}}
	\def\bl{\begin{lemma}\label}
		\def\el{\end{lemma}}
	\def\bp{\begin{proposition}\label}
		\def\ep{\end{proposition}}
	\def\bex{\begin{example}\label}
		\def\ex{\end{example}}
	\def\bas{\begin{assumption}}
		\def\eas{\end{assumption}}
	\def\br{\begin{remark}\label}
		\def\er{\end{remark}}
	\def\ba{\begin{array}}
		\def\ea{\end{array}}
	\def\ed{\end{document}}

\def\rf{\eqref}

\def\square#1{\vbox{\hrule\hbox{\vrule height#1%
			\kern#1\vrule}\hrule}}
\def\rectangle#1#2{\vbox{\hrule\hbox{\vrule height#1%
			\kern#2\vrule}\hrule}}

\font\tenbb=msbm10 \font\sevenbb=msbm7 \font\fivebb=msbm5

\newfam\bbfam
\scriptscriptfont\bbfam=\fivebb \textfont\bbfam=\tenbb
\scriptfont\bbfam=\sevenbb

\newtheorem{theorem}{Theorem}[section]
\newtheorem{corollary}[theorem]{Corollary}

\newtheorem{lemma}[theorem]{Lemma}
\newtheorem{proposition}[theorem]{Proposition}
\newtheorem{assumption}[theorem]{Assumption}

\theoremstyle{definition}
\newtheorem{definition}[theorem]{Definition}

\newtheorem{remark}[theorem]{Remark}

\newtheorem{example}{Example}[section]

\makeatletter

\@addtoreset{equation}{section}
\makeatother

\usepackage{xcolor}
\makeatletter

\input pdf-trans
\newbox\qbox
\def\usecolor#1{\csname\string\color@#1\endcsname\space}
\newcommand\bordercolor[1]{\colsplit{1}{#1}}
\newcommand\fillcolor[1]{\colsplit{0}{#1}}
\newcommand\outline[1]{\leavevmode%
	\def\maltext{#1}%
	\setbox\qbox=\hbox{\maltext}%
	\boxgs{Q q 2 Tr \thickness\space w \fillcol\space \bordercol\space}{}%
	\copy\qbox%
}
\makeatother
\newcommand\colsplit[2]{\colorlet{tmpcolor}{#2}\edef\tmp{\usecolor{tmpcolor}}%
	\def\tmpB{}\expandafter\colsplithelp\tmp\relax%
	\ifnum0=#1\relax\edef\fillcol{\tmpB}\else\edef\bordercol{\tmpC}\fi}
\def\colsplithelp#1#2 #3\relax{%
	\edef\tmpB{\tmpB#1#2 }%
	\ifnum `#1>`9\relax\def\tmpC{#3}\else\colsplithelp#3\relax\fi
}
\bordercolor{black}
\fillcolor{white}
\def\thickness{.3}

\def\hTh{\outline{$\Theta$}}
\title{Long-Time Behaviors of Stochastic Linear-Quadratic Optimal Control Problems}

\author{Jiamin Jian
\thanks{Department of Mathematics, University of Michigan, Ann Arbor, MI 48109, USA. {\tt jiaminj@umich.edu.}}
~~~~Sixian Jin
\thanks{Department of Mathematics,
California State University San Marcos, San Marcos, CA 92096,
USA. {\tt sjin@csusm.edu.}}
~~~~Qingshuo Song
\thanks{Department of Mathematical Sciences,
Worcester Polytechnic Institute, Worcester, MA 01609, USA. {\tt qsong@wpi.edu.}}
~~~~Jiongmin Yong \thanks{Department of Mathematics, University of Central Florida, Orlando, FL 32816, USA, {\tt jiongmin.yong@ucf.edu}.
                    This author is supported in part by NSF grant DMS-2305475. }
}

\date{}

\begin{document}

\maketitle


\begin{abstract}
This paper investigates the asymptotic behavior
of the solution to a linear-quadratic stochastic optimal control problems. 
The so-called probability cell problem is introduced the first time. 
It serves as the probability interpretation of the  well-known cell problem 
in the homogenization of Hamilton-Jacobi-Bellman equations. By establishing a connection 
between this problem and the ergodic cost problem, we reveal the turnpike 
properties of the linear-quadratic stochastic optimal control problems 
from various perspectives.
\end{abstract}


\section{Introduction}
\label{s:introduction}

One interesting feature of (deterministic) 
optimal control problem in large time duration 
is that under certain conditions, optimal pair will be close to that of some (static) optimization problem in some sense. This is called the {\it turnpike property} of the optimal control problem, suggested by the highway system in the USA. Ramsey \cite{Ramsey 1928} was the first who realized such
a fact when he was studying economic growth problems in an infinite horizon in 1928. In 1945, von Neumann \cite{Neumann 1945} further developed the relevant theory. The name {\it turnpike property} was coined by Dorfman--Samuelson--Solow in 1958 (see \cite{Dorfman-Samuelson-Solow 1958}). Since then, the turnpike phenomenon has gained significant attention for finite and infinite dimensional problems in the context of deterministic discrete-time and continuous-time systems, see relevant papers \cite{McKenzie 1976, Zaslavski 2011, Porretta-Zuazua 2013, Trelat-Zuazua 2015, Porretta-Zuazua 2016, Lou-Wang 2019, Breiten-Pfeiffer 2020, Esteve-Kouhkouh-Pighin-Zuazua 2020}, and books \cite{Carlson-Haurie-Leizariwitz 1991, Zaslavski 2005, Zaslavski 2019} for nice updated surveys. The recent paper \cite{Sun-Wang-Yong 2022} by Sun--Wang--Yong in 2022 was the first that addresses the turnpike property for stochastic linear-quadratic (LQ) optimal control problems. 
Moreover, in \cite{Sun-Yong 2024}, Sun--Yong further established the exponential, integral, and mean-square turnpike properties for optimal pairs of linear quadratic mean-field stochastic control problems, under the stabilizability condition for the homogeneous state equation.

\ms

To be more precise, let us consider a one-dimensional standard Brownian motion $\{W(t)\}_{t \ges0}$ defined on a complete filtered probability space $(\O,\cF,\dbF,\dbP)$ which satisfies the usual conditions, with $\dbF=\{\cF_t\}_{t\ges0}$ being the natural filtration of $W(\cd)$ augmented by all the $\dbP$-null sets in $\cF$. Let $T>0$ be a large time horizon. Consider the following controlled linear stochastic differential equation (SDE) on time horizon $[t_0,T]$:
\bel{state}\left\{\2n\ba{ll}
\ns\ds dX(t)=\big[AX(t)+Bu(t)+b\big]dt+\big[CX(t)+Du(t)+\si\big]dW(t),\qq t\in[t_0,T],\\
\ns\ds X(t_0)=x,\ea\right.\ee
where $A,C\in\dbR^{n\times n}$, $B,D\in\dbR^{n\times m}$, $b,\si,x\in\dbR^n$, with $\dbR^{n\times m}$ being the set of all $n\times m$ real matrices, and $\dbR^n$ being the standard $n$-dimensional real Euclidean space. Clearly, for any initial pair $(t_0,x)\in[0,T]\times\dbR^n$ and
\bel{U[s,T]}\ba{ll}
\ns\ds u(\cd)\in\sU[t_0,T]\equiv\Big\{u:[t_0,T]\times\O\to\dbR^m\bigm|u(\cd)\hb{ is $\dbF$-progressively measurable},\\
\ns\ds\qq\qq\qq\qq\qq\qq\qq\qq\qq\dbE \[\int_{t_0}^T|u(t)|^2dt \]<\i\Big\},\ea\ee
the state equation \rf{state} admits a unique solution
$$\ba{ll}
\ns\ds X(\cd)=X(\cd\,;t_0,x,u(\cd))\in
L^2_\dbF(t_0,T;\dbR^n)\\
\ns\ds\qq\;\equiv\Big\{\f:[t_0,T]\times\O\to\dbR^n\bigm|\f(\cd)\hb{ is $\dbF$-progressively measurable with continuous paths,}\\
\ns\ds\qq\qq\qq\qq\qq\qq\qq\dbE\[\sup_{t\in[t_0,T]}|\f(t)|^2\]
<\i\Big\}.\ea$$
In what follows, let $\dbS^n$ ($\dbS^n_+$, $\dbS^n_{++}$)
be the collection of
real symmetric (positive semi-definite, positive definite) matrices. To measure the performance of the control, we introduce the following quadratic cost functional:
\bel{cost}J^T(t_0,x;u(\cd))=\dbE \[\int_{t_0}^Tf(X(t),u(t))dt\],\ee
where
\bel{f}f(x,u)=\frac{1}{2}\big(\lan Qx,x\ran+2\lan Sx,u\ran+\lan Ru,u\ran+2\lan q,x\ran+2\lan r,u\ran\big),\ee
with $Q\in\dbS^n_{++}$, $R\in\dbS^m_{++}$, $S\in\dbR^{m\times n}$, satisfying $Q-S^\top R^{-1}S\in\dbS^n_{++}$, and $q\in\dbR^n$, $r\in\dbR^m$. The {\it finite-horizon stochastic LQ optimal control problem} can be formulated as follows.

\ms

{\bf Problem (LQ)$^T$.} For any initial pair $(t_0,x)\in[0,T]\times\dbR^n$, find an {\it open-loop optimal control} $u^T(\cd)\in\sU[t_0,T]$ such that
\bel{eq:value_function}V^{T}(t_0,x)\equiv J^T \big(t_0,x;u^T(\cd) \big)=\inf_{u(\cd)\in\sU[t_0,T]}J^T(t_0,x;u(\cd)).\ee
The corresponding {\it open-loop optimal state process}
is denoted by $X^T(\cd)\equiv X(\cd\,;t_0,x,u^T(\cd))$.
We call $V^T(\cd\,,\cd)$ the {\it value function} of Problem (LQ)$^T$. Under the above-mentioned assumptions, functional $u(\cd)\mapsto J^T(t_0,x;u(\cd))$ is uniformly convex, and Problem (LQ)$^T$ admits a unique {\it optimal pair} $(X^T(\cd),u^T(\cd))$. When this happens, Problem (LQ)$^T$ is said to be ({\it uniquely}) {\it open-loop solvable} at $(t_0,x)\in[0,T]\times\dbR^n$. In the case $t_0 = 0$, we use the notation $V^T(x) \equiv V^T(0, x)$ for simplicity.

\ms

When $C, D$ and $\si$ are all zero with suitable dimensions, the {\it deterministic turnpike property} is well-known, and some interesting results were established in \cite{Porretta-Zuazua 2013}. It says that (with the initial time $t_0=0$) there exist some absolute constants $\l,K>0$, and $(x^*,u^*)\in\dbR^n\times\dbR^m$
such that
\bel{eq:tp01} |X^T(t)-x^*|+|u^T(t)-u^*|\les K \big(e^{-\l t}+e^{-\l(T-t)} \big), \qq  \forall t\in[0,T],\ee
where $(x^*,u^*)$ is the solution of the following {\it static optimization problem}:
\bel{D-static}\left\{\2n\ba{ll}
\ds\hb{minimize }L(x,u)\equiv \frac{1}{2}\big(\lan Qx,x\ran+2\lan Sx,u\ran+\lan Ru,u\ran+2\lan q,x\ran+2\lan r,u\ran\big)\\
\ns\ds\hb{subject to }Ax+Bu+b=0,\ea\right.\ee
In the above, $|\,\cd\,|$ stands for the Euclidean norm of $\dbR^n$ or $\dbR^m$. This implies
\bel{<e^T}|X^T(t)-x^*|+|u^T(t)-u^*|\les2Ke^{-\l\d T},\qq \forall t\in[\d T, (1-\d)T],\ \d\in(0,1/2).\ee
When $T$ is large, the above estimate tells that the open-loop optimal pair $(X^T(\cd),u^T(\cd))$ can be well approximated by $(x^*,u^*)$ for $t\in[\d T,(1-\d)T]$. For example, we may take $\d=\frac{1}{100}$. Then in the interval $[\frac{T}{100}, \frac{99}{100}T]$, a major part of $[0,T]$, estimate \rf{<e^T} holds with the right-hand side $2Ke^{-\frac{\l}{100}T}$ being very small.

\ms

In the recent work \cite{Sun-Yong 2024},  
Sun--Yong have extended the turnpike property to the stochastic LQ optimal control problem, where at least one of $C$, $D$, and $\si$ is nonzero. This extension is by no means trivial and poses challenges because the presence of Brownian noise makes it impossible for any control to freeze the state at a fixed point. 
It was proved that (see Theorem 3.2 of \cite{Sun-Yong 2024}) (with $t_0=0$ again) for some absolute constants $K,\l>0$:
\bel{stoch-turnpike}\dbE\big[|X^T(t)-\BX^*(t)|^2+|u^T(t)-\Bu^*(t)|^2\big]\les K\big(e^{-\l t}+e^{-\l(T-t)} \big),\qq \forall t\in[0,T],\ee
where $\BX^*(\cd)$ and $\Bu^*(\cd)$ are two constructed stochastic processes independent of $T$, satisfying
\bel{E}\dbE[\BX^*(t)]=x^*,\qq\dbE[\Bu^*(t)]=u^*,\ee
with $(x^*,u^*)$ being the solution of the following {\it static optimization problem}:
\bel{S-static}\left\{\2n\ba{ll}
\ds\hb{minimize }L(x,u)\equiv \frac{1}{2}\big(\lan Qx,x\ran+2\lan Sx,u\ran+\lan Ru,u\ran+2\lan q,x\ran+2\lan r,u\ran\big)\\
\ns\ds\qq\qq\qq\qq\qq+\frac{1}{2}\lan P(Cx+Du+\si),Cx+Du+\si\ran,\\
\ns\ds\hb{subject to }Ax+Bu+b=0.\ea\right.\ee
In the above, $P\in\dbS^n_{++}$ is a solution
to the algebraic Riccati equation
\bel{ARE}\ba{ll}
\ns\ds PA+A^\top P+C^\top PC+Q\\
\ns\ds \qq\qq-(PB+C^\top PD+S^\top)(R+D^\top PD)^{-1}(PB+C^\top PD+S^\top)^\top=0,\ea\ee
satisfying certain conditions. It is worth noting that in the deterministic case, i.e., $C=0$, $D=0$, and $\si=0$, \rf{S-static} is reduced to \rf{D-static}. The above estimate \rf{stoch-turnpike} is referred to as {\it stochastic turnpike property}.

\ms

\ms

For Problem (LQ)$^T$, we know that the value function $V^T(\cd\,,\cd)$ is of the following form:
\bel{V^T} V^T(t,x)=\lan P^T(t)x,x\ran+2\lan p^T(t),x\ran+p_0^T(t),\qq\forall(t,x)\in[0,T]\times\dbR^n, \ee
for some differentiable functions $P^T(\cd)$, $p^T(\cd)$, and $p_0^T(\cd)$.
Next, let us define Hamiltonian as follows:
\bel{H}\ba{ll}
\ns\ds\dbH(x,\Bp,\BP,u)=\lan Ax+Bu+b,\Bp\ran+\frac{1}{2}\lan\BP(Cx+Du+\si),Cx+Du+\si\ran\\
\ns\ds\qq\qq\qq\qq\qq+\frac{1}{2}\big(\lan Qx,x\ran+2\lan Sx,u\ran+\lan Ru,u\ran+2\lan q,x\ran+2\lan r,u\ran\big),\\
\ns\ds H(x,\Bp,\BP)=\inf_{u\in\dbR^m}\dbH(x,\Bp,\BP,u).\ea\ee
Then by dynamic programming principle, since the value function $V^T(\cd\,,\cd)$ is smooth, it is a classical solution of the following Hamilton-Jacobi-Bellman (HJB) equation:
\bel{HJB}\left\{\2n\ba{ll}
\ns\ds V^T_t(t,x)+H \big(x,V^T_x(t,x),V^T_{xx}(t,x) \big) = 0, \q (t,x)\in[0,T]\times\dbR^n,\\
\ns\ds V^T(T,x)=0,\q x\in\dbR^n,\ea\right.\ee
where $V^T_x(t,\cd)$ and $V^T_{xx}(t,\cd)$
are gradient and Hessian of $V^T(t,\cd)$ in $x$, respectively.

\ms

Recall that the classical Feynman-Kac formula says that the solution $V^T(\cd\,,\cd)$ of the PDE (like \rf{HJB}) can be identified as the value function of Problem (LQ)$^T$. Roughly speaking, Problem (LQ)$^T$ gives a probability interpretation of the HJB equation \rf{HJB}.

\ms

Now, we consider the infinite-horizon case with $b=\si=0$ in \rf{state} and $t_0 = 0$, i.e., we consider the homogeneous state equation
\bel{state2}\left\{\2n\ba{ll}
\ns\ds dX(t)=\big[AX(t)+Bu(t)\big]dt+\big[CX(t)+Du(t)\big]dW(t),\q t\ges0,\\
\ns\ds X(0)=x,\ea\right.\ee
and the cost functional
$$J^\i_0(x;u(\cd))=\dbE\[\int_0^\i f^0(X(t),u(t))dt\],$$
where ($q=0$ and $r=0$ in \rf{f})
$$f^0(x,u)=\frac{1}{2}\big(\lan Qx,x\ran+2\lan Sx,u\ran+\lan Ru,u\ran\big).$$
Correspondingly, we define (compare to \rf{U[s,T]})
\bel{U[s,i]}\ba{ll}
\ns\ds\sU[0,\i)\equiv\Big\{u:[0,\infty)
\times\O\to\dbR^m\bigm|u(\cd)
\hb{ is $\dbF$-progressively measurable},\\
\ns\ds\qq\qq\qq\qq\qq\qq\qq\qq\qq\dbE\[\int_0^\i|u(t)|^2dt\]<\i\Big\},\ea\ee
and
 \bel{U_ad[0,i)}\sU_{ad}[0,\i)\equiv\big\{u(\cd)\in\sU[0,\i)\bigm|X(\cd\,;x,u(\cd))\in L^2_\dbF(0,\i;\dbR^n)\big\},
\ee
where
$$\ba{ll}
L^2_\dbF(0,\infty;\dbR^n)\equiv\Big\{\f:[0,\infty)\times\O\to\dbR^n\bigm|\f(\cd)\hb{ is $\dbF$-progressively measurable with continuous paths,}\\
\ns\ds\qq\qq\qq\qq\qq\qq\qq\qq\qq\dbE\[\sup_{t\in[0,\infty)}|\f(t)|^2\]
<\i\Big\}.\ea$$
Clearly, in general, $\sU_{ad}[0,\i)\ne\sU[0,\i)$, and it might even be that $\sU_{ad}[0,\i)=\varnothing$. When system \rf{state2}, denoted by $[A,C;B,D]$, is stabilizable, one can show easily that $\sU_{ad}[0,\i)\ne\varnothing$ (see below). Hence, in this case we can formulate the following {\it homogeneous LQ control problem}.

\ms

{\bf Problem (LQ)$^\i_0$.} For any $x\in\dbR^n$, find an {\it open-loop optimal control} $u^\i(\cd)\in\sU_{ad}[0,\i)$ such that
\bel{V^i}V^\i(x)=J_0^\i(x;u^\i(\cd))=\inf_{u(\cd)\in\sU_{ad}[0,\i)}J_0^\i(x;u(\cd)).\ee

Similar to Problem (LQ)$^T$,
we denote the corresponding
{\it open-loop optimal state process}
by $X^\i(\cd)\equiv X(\cd\,;x,u^\i(\cd))$,
 and call $V^\i(\cd)$ the value function of
 Problem (LQ)$^\i_0$, $(X^\i(\cd), u^\i(\cd))$ the
{\it optimal pair} of the problem, respectively.
Similar to the finite-horizon case, we have
\bel{V^i2}V^\i(x)=\lan Px,x\ran,\qq\forall x\in\dbR^n,\ee
for some $P\in\dbS^n$. Next, let the corresponding Hamiltonian be defined as follows:
\bel{H^0}\ba{ll}
\ns\ds\dbH^0(x,\Bp,\BP,u)=\lan Ax+Bu,\Bp\ran+\frac{1}{2}\lan\BP(Cx+Du),Cx+Du\ran\\
\ns\ds\qq\qq\qq\qq\qq+\frac{1}{2}\big(\lan Qx,x\ran+2\lan Sx,u\ran+\lan Ru,u\ran\big),\\
\ns\ds H^0(x,\Bp,\BP)=\inf_{u\in\dbR^m}\dbH^0(x,\Bp,\BP,u).\ea\ee
Then it follows from the dynamic programming principle that since the value function $V^\i(\cd)$ is smooth, it is a classical solution to the following HJB equation:
\bel{HJB2}H^0(x,V^\i_x(x),V^\i_{xx}(x))=0,\qq x\in\dbR^n,\ee
where $V^\i_x(\cd)$ and $V^\i_{xx}(\cd)$ are gradient and Hessian of $V^\i(\cd)$, respectively. Thus, from the classical Feynman-Kac formula viewpoint again, roughly speaking, Problem (LQ)$^\i_0$ is a probability interpretation of HJB equation \rf{HJB2}.

\ms

Now, if we allow one of $b,\si,q,r$ to be nonzero, then we may solve the 
nonhomogeneous state equation \rf{state} on $[0,\i)$. But, the cost functional
$$J^\i(x;u(\cd))=\dbE\[\int_0^\i f(X(t),u(t))dt\]$$
might not be convergent, regardless whether the homogeneous system is stabilizable. 
Consequently, the corresponding LQ control problem, denoted formally by (LQ)$^\i$, is not well-formulated. However, the so-called {\it cell problem} under the LQ framework introduced in \cite{Lions-Papanicolauo-Varadhan 1987}, which is closely related to the homogenization of HJB equations, is always meaningful. In fact, by definition, the cell problem can be formulated as follows.

\ms

\bf Problem (C). \rm For a Hamiltonian $H(\cd, \cd, \cd)$ defined in \eqref{H}, seek a pair $(V(\cd), c_0)\in C^2(\dbR^n)\times\dbR$ such that for all $x\in\dbR^n$,
\bel{cell}H\big(x, V_x(x), V_{xx}(x)\big)=c_0.\ee

\rm

\bs

At a first look of Problem (C), one can formally plug any $V(\cd) \in C^2(\dbR^n)$ 
into the Hamiltonian to obtain a corresponding $c_0$. However, the key requirement 
is that the resulting expression $c_0$ must be independent of $x$. 
Then, a natural question arises: 
What is the probability interpretation of the above Problem (C) in the connection 
to the ill-posed stochastic LQ control problem under infinite horizon with non-zero coefficients?

\ms

Clearly, it cannot be Problem (LQ)$_0^\i$ as $b,\si,q,r$ are allowed to be non-zero. It is not Problem (LQ)$^\i$ either, as this problem is not even well-formulated.
Also, since the problem seems to be relevant to LQ control problem posed in infinite horizon, we expect that the problem should be closely related to long time behavior of optimal pair to the corresponding LQ control problem. Consequently, it is expected that the turnpike property of stochastic LQ control problem will somehow plays an interesting role here.

\ms
The purpose of the current paper is mainly to investigate the above question and beyond. 
The main contributions can be summarized as follows:

\ms

$\bullet$ For HJB equations \rf{HJB} and \rf{HJB2} as well as the cell problem \rf{cell}, resulting from LQ control problems, classical solutions will be calculated. It turns out that all these solutions are up to quadratic in $x$, involving differential/algebraic Riccati equations, and terminal value problems for ODEs, whose solvability follows under proper conditions.

\ms

$\bullet$ According to the results in \cite{Sun-Yong 2020}, for Problem (LQ)$^T$, under proper conditions, the open-loop optimal state process $X^T(\cd)$ can be determined by the closed-loop system, and the open-loop optimal control $u^T(\cd)$ admits a closed-loop representation, involving differential Riccati equation and a terminal value problem of ODE. The coefficients in these equations have natural convergence as $T\to\i$. The crucial step to show these convergence is to prove the uniform boundedness of the solution $p^T(\cd)$ of the involved ODE, see Theorem \ref{approximation}.

\ms

$\bullet$ The so-called probabilistic cell problem is introduced the first time. It turns out that we are the first to bridge the connection between the cell problem (Problem (C)) originated from the PDE literature and a proper stochastic control problem (Problem (PC)). In particular, our nonconventional formulation of Problem (PC) leads to a new verification theorem that provides a stochastic representation for the solution to the cell problem, see Theorem \ref{t:verify}.


\ms

$\bullet$ The probabilistic cell problem provides an alternative
method for establishing the stochastic turnpike property \rf{stoch-turnpike} (see Theorem \ref{t:convergence_of_optimal_path} and compare with \cite{Sun-Yong 2024}). Moreover, we provide the connection between the probabilistic cell problem and
the {\it ergodic cost problem}
(see Remark 3.6.7 in \cite{Arapostathis-Borkar-Ghosh 2012}),
which involves determining the constant
\footnote{From now on, $\bar c$ will be
used to denote the ergodic constant in \eqref{eq:ecost1},
while $c_0$ is a part of the solution to the cell problem \eqref{cell}.
They are not necessarily the same from the definition and
it is one of our objective to prove the connection between them.}
\bel{eq:ecost1}
\bar c
\equiv\lim_{T\to\i}\frac{1}{T}V^T(x)=\lim_{T\to\i}\frac{1}{T}J^T(x; u^T(\cd)).\ee
Therefore, stochastic ergodic optimal control problem becomes another probability interpretation of the cell problem. This unveils a new turnpike property in terms of the cost function in addition to the aforementioned turnpike property of \rf{stoch-turnpike} with respect to the control process and state process:
\bel{eq:tpk03}\lim_{T\to\i}\frac{1}{T}J^T(x; \bar u(\cd))=\bar c,\ee
where $\bar u(\cd)$ is a $T$-independent control process obtained from the probabilistic cell problem, see Theorem \ref{t:convergence_of_value_function}.
\ms

$\bullet$ We fill in the connections between the stochastic turnpike properties and the static optimization problem in Corollary \ref{c:connection} by explicitly writing out the relationship between the points $(x^*,u^*)$ and the coefficients of the optimal control from the classical cell problem.
\ms

$\bullet$ In addition, our framework establishes the turnpike property for arbitrary initial states in $\mathbb R^n$. As a byproduct, this new framework can recover the results in \cite{Sun-Wang-Yong 2022, Sun-Yong 2024} by making a special choice for the initial state of $\bar{X}(\cdot)$.
\ms

The five problems, namely, the {\it finite time stochastic control problem} \rf{eq:value_function}, the {\it static optimization problem} \rf{S-static}, the {\it cell problem} \rf{cell}, the {\it probabilistic cell problem} \rf{eq:ergo_lq}, and the {\it ergodic cost problem} \rf{eq:ecost1} will be interwoven throughout the remainder of this paper in the following manner:
In Section \ref{s:LQ}, we recall the existing results on the finite- and infinite-horizon stochastic LQ control problems, including the open-loop and closed-loop solvability under the two main assumptions (H1)--(H2).
In Section \ref{s:cell_problem}, we discuss the solutions to the HJB equations of the finite- and infinite-horizon LQ control problems, and then introduce the cell problem, solving for its solution under our main assumptions (H1)--(H2).
In Section \ref{s:NE}, we compare the parameters obtained from the previous
two sections and establish several interesting convergence results,
which play a crucial role for the establishment of turnpike properties.
These estimates also introduce a meaningful connection between the cell problem and the static optimization problem. Section \ref{s:probabilistic_cell_problem} introduces
the probabilistic cell problem and proves its verification theorem.
Section \ref{s:convergence} presents the new turnpike properties of the optimal pairs and the cost functions between the finite-horizon stochastic control problem and the probabilistic cell problem.

\section{Finite- and Infinite-Horizon Stochastic LQ Control Problems}
\label{s:LQ}
In this section, we recall the finite- and infinite-horizon
stochastic LQ control problems. Most of the results
come from the prior work of one of the authors,
and proofs can be found in references
\cite{Huang-Li-Yong 2015, Sun-Yong 2020, Sun-Wang-Yong 2022, Sun-Yong 2024}.


\subsection{Finite-horizon LQ control problem: (LQ)$^T$}

\ms

Let us recall the following. In the previous section, we have formulated Problem (LQ)$^T$. 
Next, let $\hTh[t_0,T]=L^\i(t_0,T;\dbR^{m\times n})$. For any $(\Th(\cd),\theta(\cd))\in\hTh[t_0,T]\times\sU[t_0,T]$, let
$$u(t)=\Th(t)X(t)+\theta(t), \qq t\in[t_0,T],$$
with $X(\cd) \equiv X(\cd\,;t_0,x,\Th(\cd),\theta(\cd))$ being the solution to the following {\it closed-loop} system
\bel{closed}\left\{\2n\ba{ll}
\ds dX(t)=\big[\big(A+B\Th(t)\big)X(t)+B\th(t)+b\big]dt\\
\ns\ds\qq\qq\qq+\big[\big(C+D\Th(t)\big)X(t)+D\th(t)+\si\big]dW(t),\q t\in[t_0,T],\\
\ns\ds X(t_0)=x.\ea\right.\ee
We adopt the following notation:
$$J^T(t_0,x;\Th(\cd),\theta(\cd))=J^T(t_0,x;\Th(\cd)X(\cd)+\theta(\cd)).$$
A couple $(\Th^T(\cd),\th^T(\cd))$ is called a {\it closed-loop optimal strategy} at $t_0$, if
$$J^T(t_0,x;\Th^T(\cd),\th^T(\cd))\les J^T(t_0,x;\Th(\cd),\theta(\cd)),\qq \forall (\Th(\cd),
\theta(\cd))\in\hTh[t_0,T]\times\sU[t_0,T],~x\in\dbR^n.$$
When such a pair $(\Th^T(\cd),\th^T(\cd))$ exists, Problem (LQ)$^T$ is said to be {\it closed-loop solvable} on $[t_0,T]$.

\ms

Next, we introduce the following hypothesis.

\ms

{\bf(H1)} The matrices $Q\in\dbS^n_{++}$ and $R\in\dbS^m_{++}$ with $Q-S^\top R^{-1}S\in\dbS^n_{++}$.

\ms

To proceed further, we present the following differential Riccati equation
\bel{DRE}\left\{\2n\ba{ll}
\ns\ds\dot P^T(t)+P^T(t)A+A^\top P^T(t)+C^\top P^T(t)C+Q-[P^T(t)B+C^\top P^T(t)D+S^\top ]\\
\ns\ds\qq\qq\cd[R+D^\top P^T(t)D ]^{-1} [B^\top P^T(t)+D^\top
P^T(t)C+S]=0, \q t\in[0,T],\\
\ns\ds P^T(T)=0,\ea\right.\ee
and the terminal value problem for a backward ODE (BODE, for short)
\bel{p^T}\left\{\2n\ba{ll}
\ns\ds\dot p^T(t)+[A+B\Th^T(t)]^\top p^T(t)+[C+D\Th^T(t)]^\top P^T(t)\si+P^T(t)b\\
\ns\ds\qq\qq\qq+q+\Th^T(t)^\top r=0,\q t\in[0,T],\\
\ns\ds p^T(T)=0.\ea\right.\ee
For Problem (LQ)$^T$, we recall the following known results; see Theorem 2.3.3 and Corollary 2.5.7 in \cite{Sun-Yong 2020} for the open-loop and closed-loop solvability, respectively.

\ms

\bl{l:finite_time_control} \sl Let {\rm(H1)} hold. Then

\ms

{\rm(i)} The pair $(X^T(\cd),u^T(\cd))$ is the unique open-loop optimal pair of Problem (LQ)$^T$ at every $(t_0,x)\in[0,T)\times\dbR^n$ if and only if the following forward-backward SDE (FBSDE, for short) admits adapted solution $(X^T(\cd),Y^T(\cd),Z^T(\cd))$:
\bel{FBSDE}\left\{\2n\ba{ll}
\ns\ds dX^T(t)=[AX^T(t)+Bu^T(t)+b]dt+[CX^T(t)+Du^T(t)+\si]dW(t),\\
\ns\ds dY^T(t)=-[A^\top Y^T(t)+C^\top Z^T(t)+QX^T(t)+S^\top u^T(t)+q]dt+Z^T(t)dW(t),\\
\ns\ds X^T(t_0)=x, \,\, Y^T(T)=0,\ea\right.\ee
with stationarity condition:
\bel{stationarity}B^\top Y^T(t)+D^\top Z^T(t)+SX^T(t)+Ru^T(t)+r=0,\qq t\in[t_0,T].\ee
Moreover, the unique open-loop optimal control $u^T(\cd)$ admits the following closed-loop representation:
\bel{u^T}u^T(t)=\Th^T(t)X^T(t)+\th^T(t),\qq t\in[t_0,T],\ee
where $X^T(\cd)\equiv X(\cd\,;t_0,x,\Th^T(\cd),\th^T(\cd))$ is the solution to the corresponding closed-loop system (similar to \rf{closed}), with
\bel{Th^T}\left\{\2n\ba{ll}
\ds\Th^T(t)=-[R+D^\top P^T(t)D]^{-1}[B^\top P^T(t)+D^\top P^T(t)C+S],\\
\ns\ds\th^T(t)=-[R+D^\top P^T(t)D]^{-1}[B^\top p^T(t)+D^\top P^T(t)\si+r].\ea\right.\ee
In the above, $P^T(\cd)$ is the unique solution to the differential Riccati equation \eqref{DRE} having the property that
\bel{>d}R+D^\top P^T(t)D \ges \d I_m,\qq \forall t\in[t_0, T],\ee
for some uniform constant $\d>0$, where $I_m$ is the $m \times m$ identity matrix. Moreover, $p^T(\cd)$ is the solution to the BODE \rf{p^T}.

\ms

{\rm(ii)} The value function is given by
\bel{value[0,T]}V^T(t,x)=\frac{1}{2}\big[\lan P^T(t)x,x\ran+2\lan p^T(t),x\ran+p_0^T(t)\big],\qq (t,x)\in[t_0,T]\times\dbR^n,\ee
with
\bel{p^T_0}\ba{ll}
\ns\ds p^T_0(t)=\int_t^T\big(\lan P^T(s)\si,\si\ran+2\lan p^T(s),b\ran\\
\ns\ds\qq\q-\lan[R+D^\top\1n P^T(t)D]^{-1}[B^\top\1n p^T(s)+D^\top\1n P^T(s)\si+r],
B^\top\1n p^T(s)+D^\top\1n P^T(s)\si+r\ran\big)ds.\ea\ee

{\rm(iii)} Problem (LQ)$^T$ is closed-loop solvable at every $t\in[t_0,T]$ with the closed-loop optimal strategy $(\Th^T(\cd),\th^T(\cd))$ determined by \rf{Th^T}.

\el

Note that in the above, we have the following connection
\bel{Y^T}Y^T(t)=P^T(t)X^T(t)+p^T(t),\qq t\in[t_0,T].\ee

\subsection{Homogeneous LQ control problem in infinite horizon: (LQ)$_0^\i$}
Now, we recall the homogeneous LQ control problem in the infinite horizon $[0,\i)$.\footnote{For non-homogeneous case (namely, at least one of $b$ and $\si$ is a non-zero constant vector), the LQ control problem in an infinite horizon might be meaningless, since the cost functional might not be well-defined.} The state equation, denoted by $[A, C;B,D]$, reads ($b=\si=0$, compare with \rf{state}):
\bel{state[0,i)}\left\{\2n\ba{ll}
\ns\ds dX(t)=\big[AX(t)+Bu(t)\big]dt+\big[CX(t)+Du(t)\big]dW(t),\q t\in[0,\i),\\
\ns\ds X(0)=x.\ea\right.\ee
The cost functional reads (compare with \rf{cost}, and note that $q=0$ and $r=0$)
\bel{cost[0,i)}J^\i_0(x;u(\cd))=\lim_{T\to\i}\dbE\[\int_0^Tf^0(X(t),u(t))dt\],\ee
where the running cost is given by
$$f^0(x,u)=\frac{1}{2}\big(\lan Qx,x\ran+2\lan Sx,u\ran+\lan Ru,u\ran\big).$$
Note that the cost functional on $[0,\i)$ does not contain the terminal cost. Moreover, although we set $b=\si=0$, for any $x\in\dbR^n$ and $u(\cd)\in\sU[0,\i)$, the unique state $X(\cd)\equiv X(\cd\,;x,u(\cd))$ is not necessarily in $L^2_\dbF (0, \i;\dbR^n)$, and thus, the corresponding cost functional \rf{cost[0,i)} might not be well-defined. Thus, we need to first look at the following linear homogeneous uncontrolled SDE, denoted by $[A,C]$:
\bel{eq:homo_AC}dX(t)=AX(t)dt+CX(t)dW(t),\qq t\ges0.\ee
We recall the following results, which can be found in Propositions 3.5 and 3.6 of \cite{Huang-Li-Yong 2015} and in \cite{Sun 2024}.

\bl{l:exponential_stable} \sl The following statements are equivalent:

\ms

{\rm(i)} System $[A,C]$ is $L^2$-exponentially stable, i.e., for any $x\in\dbR^n$, the solution $X(\cd)\equiv X(\cd\,;x)$ to \rf{eq:homo_AC} satisfies
$$\lim_{t\to\i}e^{\l t}\dbE\big[|X(t)|^2\big]=0,$$
for some $\l>0$.

\ms

{\rm(ii)} System $[A,C]$ is $L^2$-globally integrable, i.e., for any $x\in\dbR^n$, the solution $X(\cd)\equiv X(\cd\,;x)$ to \rf{eq:homo_AC} satisfies
$$\int_0^\i\dbE\big[|X(t)|^2\big]dt<\i.$$

\ms

{\rm(iii)} System $[A,C]$ is $L^2$-asymptotically stable, i.e., for any $x\in\dbR^n$, the solution $X(\cd)\equiv X(\cd\,;x)$ to \rf{eq:homo_AC} satisfies
$$\lim_{t\to\i}\dbE\big[|X(t)|^2\big]=0.$$

\ms

{\rm(iv)} For any $\L\in\dbS^n_{++}$, the Lyapunov equation
\bel{Lyapunov}PA+A^\top P+C^\top PC=-\L\ee
admits a unique solution $P\in\dbS^n_{++}$.

\ms

\el

Consequently, we have the following corollary.

\bc{lya} \sl The following statements are equivalent:

\ms

{\rm(i)} The controlled system $[A,C;B,D]$ is $L^2$-exponentially stabilizable, (resp. $L^2$-globally integrable, $L^2$-asymptotically stable) i.e., there exists a $\Th\in\dbR^{m\times n}$, called a {\it $L^2$-exponential stabilizer} of the system, such that $[A+B\Th,C+D\Th]$ is $L^2$-exponentially stable (resp. $L^2$-globally integrable, $L^2$-asymptotically stable).

\ms

{\rm(ii)} For some $\Th\in\dbR^{m\times n}$ and for any $\L\in\dbS^n_{++}$, the following Lyapunov equation
\bel{Lyapunov2} P(A+B\Th)+(A+B\Th)^\top P+(C+D\Th)^\top P(C+D\Th)=-\L,\ee
admits a unique positive definite solution $P\in\dbS^n_{++}$.

\ec

Now, we introduce the following hypothesis.

\ms

{\bf(H2)} System $[A,C;B,D]$ is $L^2$-exponentially stabilizable.

\ms

Clearly, under (H1)--(H2), the following set is nonempty
$$\sU_{ad}[0,\i)\equiv\big\{u(\cd)\in\sU[0,\i)\bigm|X(\cd\,;x,u(\cd))\in L^2_\dbF(0,\i;\dbR^n)\big\}.$$
Then, for any $(x,u(\cd))\in\dbR^n\times\sU_{ad}[0,\i)$, the cost functional $J_0^\i(x;u(\cd))$ is well-defined. We see that unless $[A,C]$ is stable, for the controlled system $[A, C; B, D]$, $\sU[0,\i)\ne\sU_{ad}[0,\i)$. We now state the corresponding LQ control problem as follows:

\ms

\bf Problem (LQ)$^\i_0$. \rm For $x\in\dbR^n$, find $u^{\infty}(\cd)\in\sU_{ad}[0,\i)$ such that
\bel{infJ^i}J^\i_0(x;u^\i(\cd))=\inf_{u(\cd)\in\sU_{ad}[0,\i)}J^\i_0(x;u(\cd))\equiv V^\i(x).\ee

The above $u^\i(\cd)$ is called the {\it open-loop optimal control}, and the corresponding state process $X^\i(\cd)\equiv X(\cd\,;x,u^\i(\cd))$ is called the {\it open-loop state process}, $(X^\i(\cd),u^\i(\cd))$ and $V^\i(\cd)$ are called the {\it open-loop optimal pair} and the {\it value function}, respectively. When $u^\i(\cd)$ exists, Problem (LQ)$^\i_0$ is said to be {\it open-loop solvable} at $x\in\dbR^n$.

\ms

In what follows, we will denote the set of all stabilizers of $[A,C;B,D]$ by $\cS[A,C;B,D]$.
For any $\Th\in\cS[A,C;B,D]$, we introduce the following homogeneous closed-loop system:
\bel{closed-loop[0,i)}\left\{\2n\ba{ll}
\ns\ds dX(t)=(A+B\Th)X(t)dt+(C+D\Th)X(t)dW(t),\q t\in[0,\i),\\
\ns\ds X(0)=x.\ea\right.\ee
For convenience, we denote
$$J^\i_0(x;\Th)=J^\i_0(x;\Th X(\cd)),$$
with $X(\cd)$ being the solution to \rf{closed-loop[0,i)}. Any $\bar\Th\in\cS[A,C;B,D]$ satisfying the following is called a {\it closed-loop optimal strategy} of Problem (LQ)$_0^\i$:
$$J^\i_0(x;\bar\Th)\les J^\i_0(x;\Th), \qq \forall \Th\in\cS[A,C;B,D], \ x\in\dbR^n.$$
When the above $\bar\Th$ exists, Problem (LQ)$^\i_0$ is said to be {\it closed-loop solvable}.

\ms

The following results can be found in Theorem 5.1 of \cite{Huang-Li-Yong 2015}.

\bl{l:results_infinity_lq} \sl Let {\rm(H1)--(H2)} hold. Then

\ms

{\rm(i)} The algebraic Riccati equation \rf{ARE} admits a unique {\it stabilizing solution} $P$, i.e., $P\in\dbS^n_{++}$ solving \rf{ARE} such that
\bel{bar Th}\bar\Th\equiv -(R+D^\top PD)^{-1}(B^\top P+D^\top PC+S)\in \cS[A,C;B,D].\ee

{\rm(ii)} Problem (LQ)$^\i_0$ is open-loop solvable, with the open-loop optimal control admitting the following closed-loop representation:
\bel{bar u} u^\i(t)=-(R+D^\top PD)^{-1}(B^\top P+D^\top PC+S) X^\i (t) = \bar\Th X^\i(t),\qq t\in[0,\i),\ee
where $P$ is the stabilizing solution to \rf{ARE}, and $X^\i(\cd)$ is the solution to the closed-loop system \rf{closed-loop[0,i)}. Moreover, the value function is given by
$$V^\i(x)=\lan Px,x\ran,\qq x\in\dbR^n.$$

{\rm(iii)} Problem (LQ)$^\i_0$ is closed-loop solvable with the closed-loop optimal strategy $\bar\Th$ given by \rf{bar Th}.
\el


\section{HJB Equations and the Cell Problem}\label{s:cell_problem}

In this section, we look at the HJB equations corresponding to Problems (LQ)$^T$, (LQ)$^\i_0$, and beyond. We first recall $\dbH:\dbR^n\times\dbR^n\times\dbS^n\times\dbR^m\to\dbR$ and $H:\dbR^n\times\dbR^n\times\dbS^n\to\dbR$ from \rf{H}.
Note that the Hamiltonian $H(x,\Bp,\BP)$ is independent of $t\ges0$. Moreover, assuming $R+D^\top\BP D$ is positive definite, it is not hard to obtain that (see
\rf{H})
\bel{eq:Hamiltonian}\ba{ll}
\ns\ds H(x,\Bp,\BP)=\inf_{u\in\dbR^m}\dbH(x,\Bp,\BP,u)\\
\ns\ds\qq\qq\q=\lan A^\top\Bp,x\ran+\lan\Bp,b\ran+\frac{1}{2}\big(\lan(C^\top\BP C+Q)x,x\ran+2\lan x, C^\top\BP\si+q\ran+\lan\BP\si,\si\ran\\
\ns\ds\qq\qq\qq\q -|(R+D^\top\BP D)^{-\frac{1}{2}}[(D^\top\BP C+S)x+D^\top\BP\si+B^\top\Bp+ r]|^2\big),\ea\ee
where the minimum is attained at
\bel{hu} \h u(x, \Bp, \BP) = -(R\1n+\1n D^\top\BP D)^{-1} \big[(D^\top\BP C+S)x+D^\top\BP\si+B^\top\Bp+r \big].\ee
It is known that the HJB equation for Problem (LQ)$^T$ reads as \rf{HJB}.
Now, we have the following result.

\bp{finite} \sl Let {\rm(H1)} hold. Then the HJB equation \rf{HJB} admits a  classical solution, which is the value function of Problem (LQ)$^T$ and is given by
\bel{value[0,T]*}V^T(t,x)=\frac{1}{2}\big(\lan P^T(t)x,x\ran+2\lan p^T(t),x\ran+p^T_0(t)\big),\qq (t,x)\in[0,T]\times\dbR^n,\ee
where $P^T(\cd)$ is the solution to the differential Riccati equation \rf{DRE}, $p^T(\cd)$ is the solution of BODE \rf{p^T}, and $p^T_0(\cd)$ is given by \rf{p^T_0}.
\ep

\it Proof. \rm Let $V^T(t,x)$ be given by \rf{value[0,T]*}. Then,
$$V^T_x(t,x)=P^T(t)x+p^T(t),\q V^T_{xx}(t,x)=P^T(t).$$
Hence, suppressing $T$ and $t$ for simplicity, by \rf{DRE} and \rf{p^T}, we have
$$\ba{ll}
\ns\ds H(x,V^T_x(t,x),V^T_{xx}(t,x))=H(x,P^T(t)x+p^T(t),P^T(t))\\
\ns\ds=\lan A^\top(Px+p),x\ran+\lan Px+p,b\ran+\frac{1}{2}\big(\lan(C^\top PC+Q)x,x\ran+2\lan x,C^\top P\si+q\ran+\lan P\si,\si\ran\\
\ns\ds\qq-|(R+D^\top PD)^{-\frac{1}{2}}[(D^\top PC+S)x+B^\top(Px+p)+D^\top P\si+r]|^2\big)\\
\ns\ds=\frac{1}{2}\big[\big\lan\big(PA+A^\top P+C^\top PC+Q-(PB+C^\top PD+S^\top)(R+D^\top PD)^{-1}(B^\top P +D^\top PC+S)\big)x,x \big\ran\\
\ns\ds\qq+2\lan x,A^\top p+Pb+C^\top P\si+q-(PB+C^\top PD+S^\top)(R+D^\top PD)^{-1}(B^\top p+D^\top P\si+r)\ran\\
\ns\ds\qq+2\lan p,b\ran+\lan P\si,\si\ran-\lan(R+D^\top PD)^{-1}(B^\top p+D^\top P\si+r),B^\top p+D^\top P\si+r\ran\big]\\
\ns\ds=\frac{1}{2}\big[\big\lan\big(PA+A^\top P+C^\top PC+Q-(PB+C^\top PD+S^\top)(R+D^\top PD)^{-1}(B^\top P+D^\top PC+S)\big)x,x \big\ran\\
\ns\ds\qq+2\lan x,(A+B\bar\Th)^\top p+Pb+(C+D\bar\Th)^\top P\si+q+\bar\Th^\top r \ran\\
\ns\ds\qq+2\lan p,b\ran+\lan P\si,\si\ran-\lan(R+D^\top PD)^{-1}(B^\top p+D^\top P\si+r),B^\top p+D^\top P\si+r\ran\big]\\
\ns\ds=-\frac{1}{2}\big(\lan\dot P^T(t)x,x\ran+2\lan\dot p^T(t),x\ran+\dot p^T_0(t)\big)=-V_t(t,x)\ea$$
for all $(t,x)\in[0,T]\times\dbR^n$. Then, our conclusions follow.
\endpf

\ms

By taking $b=\si=q=0$ and $r=0$ in \rf{H} and \rf{eq:Hamiltonian}, we have (see \rf{H^0})
\bel{H^0_explicit}\ba{ll}
\ns\ds H^0(x,\Bp,\BP)=\inf_{u\in\dbR^m}\dbH^0(x,\Bp,\BP,u)\\
\ns\ds\qq\qq\q\;=\lan A^\top\Bp,x\ran+\frac{1}{2}\big(\lan(C^\top\BP C+Q)x,x\ran\\
\ns\ds\qq\qq\qq\qq-\big|(R+D^\top\BP D)^{-\frac{1}{2}}[(D^\top\BP C+S)x+B^\top\Bp]\big|^2\big),\ea\ee
assuming $R+D^\top\BP D$ is positive definite, then the HJB equation for Problem (LQ)$^\i_0$ reads as \rf{HJB2}. We have the following result, whose proof is a straightforward computation.

\bp{} \sl Let {\rm(H1)--(H2)} hold. Then, the HJB equation \rf{HJB2} admits a  classical solution, which is the value function $V^\i(\cd)$ of Problem (LQ)$^\i_0$ and it is given by:
$$V^\i(x)=\lan Px,x\ran,\qq x\in\dbR^n,$$
where $P$ is the stabilizing solution to the algebraic Riccati equation \rf{ARE}.

\ep

On the other hand, without assuming $b=\si=q=0$ and $r=0$, the corresponding LQ control problem, denoted by (LQ)$^\i$, is not well-formulated, since the running cost rate function $f(X(\cd),u(\cd))$ might not be integrable on $[0,\i)$, in other words, $u\notin\sU_{ad}[0,\i)$. However, the Hamiltonian $H(x,\Bp,\BP)$ is still well-defined. Hence, Problem (C) always makes sense (see \cite{Lions-Papanicolauo-Varadhan 1987}). Note that this problem is closely related to the homogenization of HJB equation. The next theorem gives the explicit solution of Problem (C).

\bt{p:chara_ergodic_cost1} \sl Let {\rm(H1)--(H2)} hold. Then Problem {\rm(C)} admits a solution $(V(\cd), c_0)$ with
\bel{v} V(x)=\frac{1}{2}\big(\lan Px,x\ran+2\lan p,x\ran+p_0\big),\ee
where $P\in\dbS^n_{++}$ is the unique stabilizing solution to the algebraic Riccati equation \rf{ARE} so that $A+B\bar\Th$ is invertible with $\bar\Th$ given by \rf{bar Th},
\bel{p}p=-[(A+B\bar\Th)^\top]^{-1}\big(Pb+(C+D\bar\Th)^\top P\si+q+\bar\Th^\top r\big),\ee
$p_0\in\dbR$ being an arbitrary constant, and
\bel{c_0}
c_0 =\frac{1}{2}\big(2\lan p,b\ran+\lan P\si,\si\ran-(B^\top p+D^\top P\si+r)^\top(R+D^\top PD)^{-1}(B^\top p+D^\top P\si+r)\big).\ee
\et

\it Proof. \rm First, by Lemma \ref{l:results_infinity_lq}, under (H1)--(H2), the algebraic Riccati equation \rf{ARE} admits a unique stabilizing solution $P\in\dbS^n_{++}$. It is clear that $A+B\bar\Th$ is invertible and $p$ is well-defined by \rf{p}.
Now, we define $V(\cd)$ by \rf{v}. Then,
\bel{V_x} V_x(x)=Px+p,\q V_{xx}(x)=P.\ee
Since $P$ solves equation \eqref{ARE}, together with \eqref{p}, we have
$$\ba{ll}
\ns\ds H(x,V_x(x),V_{xx}(x))=H(x,Px+p,P)\\
\ns\ds=\lan A^\top(Px+p),x\ran+\lan Px+p,b\ran+\frac{1}{2}\big(\lan(C^\top PC+Q)x,x\ran+2\lan x,C^\top P\si+q\ran+\lan P\si,\si\ran\\
\ns\ds\qq-|(R+D^\top PD)^{-\frac{1}{2}}[(D^\top PC+S)x+B^\top(Px+p)+D^\top P\si+r]|^2\big)\\
\ns\ds=\frac{1}{2}\big[\big\lan\big(PA+A^\top P+C^\top PC+Q\\
\ns\ds\qq-(PB+C^\top PD+S^\top)(R+D^\top PD)^{-1}(B^\top P+D^\top PC+S)\big)x,x \big\ran\\
\ns\ds\qq+2\lan x,A^\top p+Pb+C^\top P\si+q-(PB+C^\top PD+S^\top)(R+D^\top PD)^{-1}(B^\top p+D^\top P\si+r)\ran\\
\ns\ds\qq+2\lan p,b\ran+\lan P\si,\si\ran-\lan(R+D^\top PD)^{-1}(B^\top p+D^\top P\si+r),B^\top p+D^\top P\si+r\ran\big]\\
\ns\ds=\frac{1}{2}\big[\big\lan\big(PA+(A+B\bar\Theta)^\top P+(C+D\bar\Theta)^\top PC+Q+\bar\Th^\top S\big)x,x \big\ran\\
\ns\ds\qq+2\lan x,(A+B\bar\Th)^\top p+Pb+(C+D\bar\Th)^\top P\si+q+\bar\Th^\top r\ran\\
\ns\ds\qq+2\lan p,b\ran+\lan P\si,\si\ran-\lan(R+D^\top PD)^{-1}(B^\top p+D^\top P\si+r),B^\top p+D^\top P\si+r\ran\big]=c_0
\ea$$
for all $x\in\dbR^n$. Then, $(V(\cd), c_0)$ is a solution of Problem (C).
\endpf

\ms

In Theorem \ref{p:chara_ergodic_cost1}, we provide the solvability of Problem (C) and constructed an explicit expression of a solution to Problem (C) with $c_0$ being in terms of a solution to algebraic Riccati equation \rf{ARE} so that $A+B\bar\Th$ is invertible. It is worth noting that Theorem \ref{p:chara_ergodic_cost1} does not imply the uniqueness of the solution to Problem (C). Indeed, on one hand, if $(V(\cd),c_0)$ is a solution, so is $(V(\cd)+K,c_0)$ for any real constant $K$. On the other hand, even $c_0$ may not be unique. A main reason is that the solution to the algebraic Riccati equation that makes $A+B\bar\Th$ invertible might not be unique (see example
below). A deeper reason might be due to the absence of compactness in the domain of $V(\cd\,,\cd)$. This may cause the major difference from the typical cell problem studied in the existing literature. More specifically, the domain of the cell problem in the literature is torus $\dbT^n$, instead of
$\dbR^n$ in our case. When the domain is compact (like $\dbT^2$), there always exists a
maximum point of the difference of the subsolution and supersolution,
together with Ishii's lemma applying to the maximum point, the comparison principle to the cell problem can be established, which leads to the uniqueness of the solution, see for instance Theorem 4.2 of \cite{Tran 2021}.
However, in our case with the domain $\dbR^n$, the maximum point of the difference of the subsolution and supersolution may not exist due to its
non-compactness. Indeed, one shall not attempt to establish
the comparison principle in this case. In the following, we provide an example to illustrate that Problem (C) possesses multiple solutions,
consequently, $c_0$ in \rf{c_0} may not be unique.

\bex{e:nonunique_P} For simplicity, we suppose $m=n$
and let $A\in\dbS^n$,
$B=R=I$, and $C=D=S=0$, where $I$ is the $n\times n$ identity matrix. Then, one can
compute that $R+D^\top PD=I$ and
$(P,p)$ satisfies
\bel{eq:P_and_p_redu}\left\{\2n\ba{ll}
\ns\ds P^2-PA-AP-Q=0,\\
\ns\ds(P-A)p-Pb+Pr-q=0.\ea\right.\ee
The first equation of \eqref{eq:P_and_p_redu} implies that
$$(P-A)^2=A^2+Q.$$
By spectral theorem, since $A^2+Q>0$, there exists multiple
choices of invertible $\D\in\dbS^n$ such that
$$A^2+Q=\D^2.$$
Accordingly, there exist multiple solutions of \eqref{eq:P_and_p_redu}
corresponding to different choices of $\D$ in the form of
$$P=A+\D,$$
and
$$p=\D^{-1}(Pb-Pr+q)$$
only if $\D$ is invertible. Therefore, Problem (C) admits multiple
solutions in the form of \rf{v} and \rf{c_0}. In particular, $c_0$ is not unique and can be given by \rf{c_0} with $P$ replaced by $A+\D$ for any different $\D\in\dbS^n$ satisfying $A^2+Q=\D^2$.

\ms

However, there exists a unique solution pair $(P, p)$ satisfying
$P - A > 0$, which is provided by Theorem \ref{p:chara_ergodic_cost1}
under (H1)--(H2).
\ex

\br{remarktheta} \rm In the case with \rf{V_x}, the minimizer $\h u$ in \rf{hu} of the Hamiltonian $\dbH$ can be simplified as
\begin{equation*}
\begin{aligned}
\h u &= -(R+D^\top PD)^{-1} \big[(D^\top PC+S)x+D^\top P\si+B^\top(Px+p)+r \big] \\
& = -(R+D^\top PD)^{-1} \big[(B^\top P+D^\top PC+S)x+D^\top P\si+B^\top p+r \big]\\
& = \bar\Th x+\bar\th,
\end{aligned}
\end{equation*}
with $\bar\Th$ given by \rf{bar Th} and
\bel{bar th}\bar\th=-(R\1n+\1n D^\top PD)^{-1}(D^\top P\si+B^\top p+r).\ee
It is worth noting that Problem (C) is posed for the case where $b,\si,q,r$ are allowed to be non-zero, which is more general than Problem (LQ)$_0^\i$.

\er

\section{Natural Estimates and Convergence}
\label{s:NE}

This section provides important estimates
between parameter functions
$P^T(\cdot)$, $p^T(\cdot)$,
$\Theta^T(\cdot)$, and $\theta^T(\cdot)$
encountered in Problem (LQ)$^T$ from Section \ref{s:LQ}, and the constant matrices/vectors $P$, $p$, $\bar{\Theta}$, and $\bar{\theta}$ encountered in Problem (C) from Section \ref{s:cell_problem}. We expect some natural convergence among them. These estimates serve as central techniques for proving the stochastic turnpike property by means of the so-called probabilistic cell problems in Section \ref{s:probabilistic_cell_problem}, as well as establishing the connection between Problem (C) and the ergodic cost problem in Section \ref{s:convergence}. Additionally, these convergences also provide interesting representations of stable constants $(x^*,u^*)$ from static optimization problems by parameters obtained from Problem (LQ)$^T$ and Problem (C). 

In what follows, $K$ and $\lambda$ denote two generic positive constants that may vary from line to line. Moreover, by abuse of notation, for any matrix $M$, we denote $|M| = \sqrt{\text{trace}(M M^\top)}$ as the Frobenius norm of $M$, where $\text{trace}(\cd)$ is the trace operator for matrices.

\subsection{Estimates between the solutions of Problem (LQ)$^T$ and Problem (C)}

Under (H1)--(H2), Theorem 4.1 of \cite{Sun-Wang-Yong 2022} and Lemma 2.3 of \cite{Sun-Yong 2024} proved that
\bel{<K}|P^T(t)-P|\les Ke^{-\l(T-t)},\qq \forall t \in [0,T],\ee
for some $K,\l>0$, independent of $T$, where $P^T(\cd)$ is the positive semi-definite solution to the differential Riccati equation \rf{DRE}
and $P$ is the stabilizing solution of algebraic Riccati equation \rf{ARE}, respectively. From \rf{<K}, it is clear that
\bel{|P|<K}|P^T(t)|\les K,\qq \forall t \in [0,T].\ee
The purpose of this section is to establish other several natural estimates. For convenience, we recall the representations of the other terms as:
\bel{barTh^T*}\left\{\2n\ba{ll}
\ns\ds\Th^T(t)=-[R+D^\top P^T(t)D ]^{-1} [B^\top P^T(t)+D^\top P^T(t)C+S ],\\
\ns\ds\th^T(t)=-[R+D^\top P^T(t)D ]^{-1} [B^\top p^T(t)+D^\top P^T(t)\si+r ];\ea\right.\ee
\bel{Th*}\left\{\2n\ba{ll}
\ns\ds\bar\Th=-(R+D^\top PD)^{-1}(B^\top P+D^\top PC+S)\in\cS[A,C;B,D],\\
\ns\ds\bar\th=-(R+D^\top PD)^{-1}(B^\top p+D^\top P\si+r);\ea\right.\ee
\bel{p^T*}\left\{\2n\ba{ll}
\ns\ds\dot p^T(t)+[A+B\Th^T(t)]^\top p^T(t)+[C+D\Th^T(t)]^\top P^T(t)\si\\
\ns\ds\qq\qq\qq+P^T(t)b+q+\Th^T(t)^\top r=0,\q t\in[0,T],\\
\ns\ds p^T(T)=0;\ea\right.\ee
and
\bel{p*}(A+B\bar\Th)^\top p+(C+D\bar\Th)^\top P\si+Pb+q+\bar\Th^\top r=0.\ee

We now state and prove the main result of this section.

\bt{approximation}
\sl Let {\rm(H1)--(H2)} hold. Then, there exist some absolute constants $K,\l>0$, independent of $T$, such that
\bel{Th-Th}|\Th^T(t)-\bar\Th|\les Ke^{-\l(T-t)},\qq \forall t\in[0,T],\ee
\bel{p-p}|p^T(t)-p|\les Ke^{-\l(T-t)},\qq \forall t\in[0,T],\ee
\bel{th-th}|\th^T(t)-\bar\th|\les Ke^{-\l(T-t)},\qq \forall t\in[0,T].\ee
Moreover, $\Th^T(\cd)$, $\th^T(\cd)$, and $p^T(\cd)$ are uniformly bounded on $[0, T]$.

\et

\it Proof. \rm By definitions, we have
\begin{equation*}
\begin{aligned}
|\Th^T(t)-\bar\Th| & \les \big|\big([R+D^\top P^T(t)D]^{-1}-(R+D^\top PD)^{-1}\big)[B^\top P^T(t)+D^\top P^T(t)C+S]\big|\\
&\qq\q+\big|(R+D^\top PD)^{-1}\big([B^\top P^T(t)+D^\top P^T(t)C]-(B^\top P+D^\top PC)\big)\big|\\
& \les Ke^{-\l(T-t)}
\end{aligned}
\end{equation*}
for all $t \in [0,T]$. This proves \rf{Th-Th} and it also implies that $\Th^T(t)$ is uniformly bounded on $[0, T]$. Next, we show \rf{p-p}. By subtracting \rf{p*} from \rf{p^T*}, we see that it is crucial to show that $p^T(\cd)$ is uniformly bounded on $[0, T]$. We now prove this. From \rf{p^T*} and uniform boundedness of $P^T(t)$ and $\Th^T(t)$, we have
$$|p^T(t)|\les\int_t^TK(|p^T(s)|+1)ds\equiv\psi(t).$$
It follows that
$$\dot \psi (t)=-K(|p^T(t)|+1)\ges-K\psi(t)-K,$$
which leads to
$$\frac{d}{dt}\big(e^{Kt}\psi(t)\big) \ges-Ke^{Kt}.$$
Then
$$-e^{Kt}\psi(t)\ges-K\int_t^Te^{Ks}ds=-e^{KT}+e^{Kt}.$$
Thus, for some absolute constant $K>0$,
$$|p^T(t)|\les\psi(t)\les e^{K(T-t)}-1,\qq \forall t\in[0,T].$$
Hence, for any fixed $t^*\in[0,T]$, one has
\bel{4.13} |p^T(t)|\les e^{Kt^*},\qq t \in [T-t^*,T].\ee
Now, we denote
\bel{bar A}A^T(t)=[A+B\Th^T(t)]^\top,\qq \bar A=A+B\bar\Th.\ee
Then, \rf{p^T*} can be written as
$$\left\{\2n\ba{ll}
\ns\ds\dot p^T(t)+A^T(t)p^T(t)+h^T(t)=0,\q t\in[0,T],\\
\ns\ds p^T(T)=0,\ea\right.$$
with
$$h^T(t)=[C+D\Th^T(t)]^\top P^T(t)\si+P^T(t)b+q+\Th^T(t)^\top r,$$
which is uniformly bounded on $[0, T]$. By assumption (H2) and Corollary \ref{lya}(ii) and taking $\L = 3I$, there exists a positive definite $\bar P$ such that the following Lyapunov inequality holds:
$$\bar P\bar A+\bar A^\top\bar P\les-3I.$$
Consequently,
\bel{eq:esti}\ba{ll}
\ns\ds\bar PA^T(t)+A^T(t)^\top\bar P=\bar P\bar A+\bar A^\top\bar P+\bar PB[\Th^T(t)-\bar\Th]+[\Th^T(t)-\bar\Th]^\top B^\top\bar P \\
\ns\ds\les \big[-3+Ke^{-\l(T-t)} \big]I\\
\ns\ds\les-2I\les\e\bar P\les I,\q\forall t\in[0,T-t^*],\ea\ee
for some fixed $t^* > 0$, and $\e>0$, assuming $T$ to be large enough (at least $T>t^*$). Now, we observe the following for $t\in[0,T-t^*]$,
$$\ba{ll}
\ns\ds\frac{d}{dt}\lan\bar Pp^T(t),p^T(t)\ran=-\lan\bar P[A^T(t)p^T(t)+h^T(t)],p^T(t)\ran-\lan\bar Pp^T(t),A^T(t)p^T(t)+h^T(t)\ran\\
\ns\ds=-\lan[\bar PA^T(t)+A^T(t)^\top\bar P]p^T(t),p^T(t)\ran-2\lan\bar Ph^T(t),p^T(t)\ran\\
\ns\ds=\e\lan\bar Pp^T(t),p^T(t)\ran-\big(\lan[\bar PA^T(t)+A^T(t)^\top\bar P]p^T(t),p^T(t)\ran+\e\lan\bar Pp^T(t),p^T(t)\ran+2\lan\bar Ph^T(t),
p^T(t)\ran\big).\ea$$
Hence, by \rf{eq:esti}, we obtain the following estimation
$$\ba{ll}
\ns\ds\lan\bar Pp^T(t),p^T(t)\ran=e^{-\e(T-t^*-t)}\lan\bar Pp^T(T-t^*),p^T(T-t^*)\ran\\
\ns\ds+\int_t^{T-t^*}\3n e^{-\e(s-t)}\big(\lan[\bar PA^T(s)+A^T(s)^\top\bar P]p^T(s),p^T(s)\ran+\e\lan\bar Pp^T(s),p^T(s)\ran+2\lan\bar Ph^T(s),
p^T(s)\ran\big)ds\\
\ns\ds\les K|p^T(T-t^*)|^2+\int_t^{T-t^*}\3n e^{-\e(s-t)}\big(-|p^T(t)|^2+\e\lan\bar Pp^T(s),p^T(s)\ran+|\bar Ph^T(t)|^2\big)ds\\
\ns\ds\les K|p^T(T-t^*)|^2+K\int_t^{T-t^*}e^{-\e(s-t)}ds\les K(1+|p^T(T-t^*)|^2),\ea$$
with an absolute constant $K>0$. Noting the positive definiteness of $\bar P$, combining \rf{4.13}, we have the uniform boundedness of $p^T(\cd)$ on $[0,T]$.

\ms

Next, we set
\bel{bar C}\h p^{\,T}(t)=p^T(t)-p, \qq \bar C = C + D \bar \Th.\ee
Then it satisfies the following (subtracting \rf{p*} from \rf{p^T*})
\begin{equation*}
\begin{aligned}
-\frac{d}{dt}\h p^{\,T}(t) & =\bar A^\top\h p^{\,T}(t)+[B(\Th^T(t)-\bar\Th)]^\top p^{\,T}(t)+\bar C^\top [P^T(t)-P]\si\\
&\qq\q+[D(\Th^T(t)-\bar\Th)]^\top P^T(t)\si+[P^T(t)-P]b+[\Th^T(t)-\bar\Th]^\top r\\
&\equiv \bar A^\top \h p^{\,T}(t)+\h h^T(t),
\end{aligned}
\end{equation*}
with $\h p^{\,T}(T)=p^T(T)-p = -p$ and
\begin{equation*}
\begin{aligned}
\h h^{\,T}(t) & = [B(\Th^T(t)-\bar\Th)]^\top p^{\,T}(t)+\bar C^\top [P^T(t)-P]\si+ [D(\Th^T(t)-\bar\Th)]^\top P^T(t)\si \\
& \qq\qq + [P^T(t)-P]b+[\Th^T(t)-\bar\Th]^\top r.
\end{aligned}
\end{equation*}
From the above, we have
$$|\h h^{\,T}(t)|\les Ke^{-\l(T-t)},\qq  \forall t\in[0,T].$$
Since $\bar A$ is stable, that is, the matrix is Hurwitz with all eigenvalues having strictly negative real parts, we may assume
\bel{e^A} |e^{\bar A^\top t}| \les Ke^{-\frac{\l}{2}t},\qq  \forall t \ges 0 \ee
for some $\l > 0$. Thus, for all $t \in [0, T]$,
\begin{equation*}
\begin{aligned}
|\h p^{\,T}(t)| &\les|e^{\bar A^\top(T-t)}p|+\int_t^T|e^{\bar A^\top(s-t)} \h h^T(s)|ds\\
&\les Ke^{-\frac{\l}{2}(T-t)}+K\int_t^Te^{-\frac{\l}{2}(s-t)}e^{-\l(T-s)}ds
\les Ke^{-\frac{\l}{2}(T-t)},
\end{aligned}
\end{equation*}
which proves \rf{p-p}. Next, by definitions, one has
\begin{equation*}
\begin{aligned}
|\th^T(t)-\bar\th| &\les\big|\big([R+D^\top P^T(t)D]^{-1}-(R+D^\top PD)^{-1}\big)[B^\top p^T(t)+D^\top P^T(t)\si+r]\big|\\
&\qq+\big|(R+D^\top PD)^{-1}\big([B^\top p^T(t)+D^\top P^T(t)\si] - (B^\top p+D^\top P\si)\big)\big|\\
&\les Ke^{-\l(T-t)},
\end{aligned}
\end{equation*}
for all $t\in[0,T]$, which is \rf{th-th}. \endpf

\ms

\subsection{Turnpike property and static optimization}
In this subsection, we apply the convergence results in Theorem \ref{approximation} to establish a significant connection between Problem (C) and the static optimization problem that was used to build up the stochastic turnpike property 
in Sun--Yong's recent work \cite{Sun-Yong 2024}. Specifically, denote $(x^*,u^*)$ as the solution to the static optimization problem \rf{S-static}, i.e.,
$$\left\{\2n\ba{ll}
\ns\ds\hb{minimize }L(x,u)\equiv\frac{1}{2}\big(\lan Qx,x\ran+2\lan Sx,u\ran+\lan Ru,u\ran+2\lan q,x\ran+2\lan r,u\ran\big)\\
\ns\ds\qq\qq\qq\qq\qq+\frac{1}{2}\lan P(Cx+Du+\si),Cx+Du+\si\ran,\\
\ns\ds\hb{subject to }Ax+Bu+b=0,\ea\right.$$
where $P\in\dbS^n_{++}$ is the stabilizing solution to the algebraic Riccati equation \eqref{ARE}. Under (H1)--(H2), we know that the feasible set is non-empty, and the objective function is coercive and convex. Hence, the optimization problem admits a unique minimizer $(x^*,u^*)\in\dbR^n\times\dbR^m$. By Lagrange multiplier method, there exists a Lagrange multiplier $y^*\in\dbR^n$ such that the following holds
\bel{Lagrange}\left\{\2n\ba{ll}
\ns\ds(Q+C^\top PC)x^*+(S+D^\top PC)^\top u^*+A^\top y^*+q+C^\top P\si=0,\\
\ns\ds(S+D^\top PC)x^*+(R+D^\top PD)u^*+B^\top y^*+r+D^\top P\si=0,\\
\ns\ds Ax^*+Bu^*+b=0.\ea\right.\ee
Now, let $(X^T(\cd),u^T(\cd))$ be the open-loop optimal pair of Problem (LQ)$^T$. Take $\bar\Th$ as \rf{bar Th}, and let $X^*(\cd)$ solve the following nonhomogeneous infinite-horizon linear SDE:
\bel{X^*}\left\{\2n\ba{ll}
\ns\ds dX^*(t)=(A+B\bar\Th)X^*(t)dt+[(C+D\bar\Th)X^*(t)+\si^*]dW(t),\qq t\ges0,\\
\ns\ds X^*(0)=0,\ea\right.\ee
with $\si^*=Cx^*+Du^*+\si$. Clearly,
$\dbE[X^*(t)]=0$. Therefore, if we set
\bel{eq:tp_yong}\ba{ll}
\ns\ds\BX^*(t)=X^*(t)+x^*,\q\Bu^*(t)=\bar\Th X^*(t)+u^*,\\
\ns\ds\BY^*(t)=PX^*(t)+y^*,\q\BZ^*(t)=P(C+D\bar\Th)X^*(t)+P\si^*,\ea\ee
then
$$\dbE[\BX^*(t)]=x^*,\q\dbE[\Bu^*(t)]=u^*,\q\dbE[\BY^*(t)]=y^*,\q\dbE[\BZ^*(t)]=P\si^*.$$
By Theorem 3.2 and Corollary 3.5 
of \cite{Sun-Yong 2024}, we have the following exponential turnpike property: there exist constants $K,\lambda>0$, independent of $T$, such that for any $t\in[0,T]$,
\bel{Turnpike_Yong}\ba{ll}
\ns\ds\dbE\big[|X^T(t)-\BX^*(t)|^2+|u^T(t)-\Bu^*(t)|^2+|Y^T(t)-\BY^*(t)|^2+|Z^T(t)-\BZ^*(t)|^2\big]\\
\ns\ds\les K\big(e^{-\l t}+e^{-\l(T-t)}\big).\ea\ee

The convergence results in Theorem \ref{approximation} imply the following interesting connection between the static optimization problem and the cell problem.

\bc{c:connection} \sl Let {\rm(H1)--(H2)} hold. Then
\bel{yu}y^*=Px^*+p,\qq u^*=\bar\Th x^*+\bar\th.\ee
Consequently,
\bel{xuy*}\left\{\2n\ba{ll}
\ns\ds x^*=-(A+B\bar\Th)^{-1}(B\bar\th+b),\\
\ns\ds u^*=-\bar\Th[(A+B\bar\Th)^{-1}(B\bar\th+b)]+\bar\th,\\
\ns\ds y^*=-P[(A+B\bar\Th)^{-1}(B\bar\th+b)]+p.\ea\right.\ee

\ec

\it Proof. \rm By \rf{Y^T}, \eqref{<K}, \rf{Turnpike_Yong}, and Theorem \ref{approximation}, we have
\begin{equation*}
\begin{aligned}
|y^*-Px^*-p| &\les 
\mathbb E \big[ |\BY^*(t)-PX^*(t)-P[\BX^*(t)-X^*(t)]-p| \big]\\
& =
\mathbb E \big[ 
	|\BY^*(t)-Y^T(t) + P^T(t)X^T(t) + p^T(t) - PX^*(t) - P[\BX^*(t)-X^*(t)] - p|
	\big] \\
& \les
\mathbb E \big[
	|\BY^*(t)-Y^T(t)|+|P^T(t)X^T(t)-P\BX^*(t)|+|p^T(t)-p|
	\big] \\
& \les K \big(e^{-\l t}+e^{-\l(T-t)} \big),\q \forall t\in [0,T].
\end{aligned}
\end{equation*}
By letting $t = T/2$ and $T\to\i$, we see that
$$y^*=Px^*+p.$$
Similarly, by \rf{u^T}, \rf{Turnpike_Yong}, and Theorem \ref{approximation},
\begin{equation*}
\begin{aligned}
|u^*-\bar\Th x^*-\bar\th| &\les 
	\mathbb E \big[
		|\Bu^*(t)-\bar\Th X^*(t)-\bar\Th[\BX^*(t)-X^*(t)]-\bar\th| 
		\big] \\
& = 
\mathbb E \big[
	|\Bu^*(t)-u^T(t)+\Theta^T(t)X^T(t)+\theta^T(t)-\bar\Th X^*(t)-\bar\Th[\BX^*(t)-X^*(t)]-\bar\th| 
	\big ]\\
& \les
\mathbb E \big[ 
	|\Bu^*(t)-u^T(t)|+|\Th^T(t)X^T(t)-\bar\Th\BX^*(t)|+|\th^T(t)-\bar\th| 
	\big]\\
& \les K \big(e^{-\l t}+e^{-\l(T-t)} \big),\q \forall t \in [0,T].
\end{aligned}
\end{equation*}
Thus, letting $t = T/2$ and $T\to\i$ again, we have
$$u^*=\Bar\Th x^*+\bar\th.$$
These prove \rf{yu}. Then, the identity
$$0=Ax^*+Bu^*+b=(A+B\bar\Th)x^*+B\bar\th+b$$ leads to $x^*=-(A+B\bar\Th)^{-1}(B\bar\th+b)$, and the rest of \eqref{xuy*} is clear.
\endpf

\ms

\br{}\rm This corollary provides us with several interesting observations. First, all expectations of four processes $\BX^*(\cd)$, $\Bu^*(\cd)$, $\BY^*(\cd)$,
and $\BZ^*(\cd)$ in \rf{eq:tp_yong} can be represented explicitly by the coefficients of the solution to Problem (C), see \eqref{xuy*}.
Note that we can roughly connect the nonhomogeneous LQ control
problem with the cell problem, and later, we will also
establish a more convincing connection with the so-called
{\it probabilistic cell problem}.
Therefore, this result brings us a nice connection between
the turnpike property and the probabilistic cell problem.
In addition, the Lagrange multiplier $y^*$, which is also
read as the mean of the adjoint process $\BY^*(\cd)$,
shares the same value as $V_x(x^*)$ with $V(\cd)$ being
the solution to Problem (C).
\end{remark}

\ms

To conclude this section, let us look at the equation that $\BX^*(\cd)$ should satisfy. Observe the following
$$\ba{ll}
\ns\ds d\BX^*(t)=(A+B\bar\Th)(\BX^*(t)-x^*)dt+\big[(C+D\bar\Th)(\BX^*(t)-x^*)+\si^*\big]dW(t)\\
\ns\ds\qq\q=\big[(A+B\bar\Th)\BX^*(t)-(A+B\bar\Th)x^*\big]dt+\big[(C+D\bar\Th)\BX^*(t) +Du^*-D\bar\Th x^*+\si\big]dW(t)\\
\ns\ds\qq\q=\big[(A+B\bar\Th)\BX^*(t)+B\bar\th+b\big]dt+\big[(C+D\bar\Th)\BX^*(t)+D\bar \th+\si\big]dW(t).\ea$$
Hence, if we let (note \rf{bar A} and \rf{bar C})
\bel{bar A bar C}\bar A=A+B\bar\Th,\q\bar C=C+D\bar\Th,\q\bar b=B\bar\th+b,\q\bar\si=D\bar\th+\si,\ee
then $\BX^*(\cd)$ is the solution to the following:
\bel{BX-SDE}\2n\left\{\2n\ba{ll}
\ns\ds d\BX^*(t)=\big[\bar A\BX^*(t)+\bar b\big]dt+\big[\bar C\BX^*(t)+\bar\si\big]dW(t),\qq t\ges0,\\
\ns\ds\BX^*(0)=x^*.\ea\right.\ee
We will see a similar system in Section \ref{s:verificantion for LQG setting} later.

\section{Probabilistic Cell Problem}
\label{s:probabilistic_cell_problem}

We have seen that Problem (LQ)$^T$ introduced in Section \ref{s:LQ} can be regarded as a probabilistic interpretation of parabolic HJB equation \rf{HJB}, and Problem (LQ)$^\i_0$ can be regarded as a probabilistic interpretation of elliptic HJB equation \rf{HJB2}, where $b=\si=q=0$ and $r=0$. In this section, we will find the similar result for Problem (C) allowing $b,\si,q,r$ to be nonzero.

\subsection{The definition of probabilistic cell problem}
\label{s:def_probabilistic_cell_problem}

Following from Section \ref{s:NE}, under assumptions (H1)--(H2), the optimal pair $(X^T(\cd),u^T(\cd))$ of Problem (LQ)$^T$ must have uniformly bounded second moments. Then, since the map $(x,u)\mapsto f(x,u)$ has no more than quadratic growth, one has
$$|J^T(x;u(\cd))|\les KT.$$
On the other hand, since $f(X(\cd),u(\cd))$ might not be integrable on $[0,\i)$,
\bel{J^i}J^\i(x;u(\cd))=\dbE \[\int_0^\i f(X(t),u(t))dt\] \ee
might not be well-defined.
Hence, in this case, Problem (LQ)$^\i$
(namely state equation \rf{state} with cost functional \rf{J^i})
might be meaningless. However, in many economics growth problems,
it is common to assume that the running cost rate function $f(X(\cd),u(\cd))$ being (strictly) positively away from 0. Thus, instead of considering the minimization of the total cost $\ds\lim_{T\to\i}J^T(x;u(\cd))$, which might not be well-defined, it is often more meaningful to minimize
average expected cost rate
\bel{V^T/T}\limsup_{T\to\i}\frac{1}{T}J^T(x;u(\cd))=
\limsup_{T\to\i}\frac{1}{T}\dbE\[\int_0^Tf(X(t),u(t))dt\].\ee
Now, by the definition of $u^T(\cd)$, we have
$$V^T(x)=\inf_{u(\cd)\in\sU[0,T]}J^T(x;u(\cd))=J^T(x;u^T(\cd)).$$
From Theorem \ref{approximation},  as $T\to\i$, we actually have
$$\ba{ll}
\ns\ds\frac{1}{T}J^T(x;u^T(\cd))=\frac{1}{T}V^T(x)\\
\ns\ds=\frac{1}{2T} \[\lan P^T(0)x,x\ran+2\lan p^T(0),x\ran+\int_0^T\big(\lan P^T(s)\si,\si\ran+2\lan p^T(s),b\ran\\
\ns\ds\qq-\lan[R+D^\top P^T(t)D]^{-1}[B^\top p^T(s)+D^\top P^T(s)\si+r],
B^\top p^T(s)+D^\top P^T(s)\si+r\ran\big)ds\]\\
\ns\ds\to\frac{1}{2}\big(\lan P\si,\si\ran+2\lan p,b\ran-\lan[R+D^\top PD]^{-1}[B^\top p+D^\top P\si+r],B^\top p+D^\top P\si+r\ran\big) \equiv c_0,\ea$$
which is given by \eqref{c_0}
in Theorem \ref{p:chara_ergodic_cost1}.
Hence, the best one can do is to select $u(\cd)$ so that the average cost rate is close to $c_0$ (from the above). Therefore, we expect to have
\bel{J^T-cT}\limsup_{T\to\i}\big(J^T(x;u(\cd))-cT\big)=\limsup_{T\to\i}\dbE\Big[\int_0^T\big(f(X(t),
u(t))-c\big)dt\Big]<\infty,\ee
for suitable control $u(\cd)$ and constant $c\,(\ges c_0)$. Further, intuitively we have the following (singular perturbation) expansion, assuming everything is fine,
\bel{expansion}J^T(x;u(\cd))=\cJ_1^T(x;u(\cd))T+\cJ_0^T(x;u(\cd))+\cJ_{-1}^T
(x;u(\cd))T^{-1}+O(T^{-2})\ee
with
$$\cJ_i^T(x;u(\cd))=\cJ_i(x;u(\cd))+o(1),\qq\hb{as }T\to\i,~i\les1,$$
for some functionals $\cJ_i^T(x;u(\cd))$ and $\cJ_i(x;u(\cd))$. Then \rf{J^T-cT} is roughly equal to $\cJ^T_0(x;u(\cd))$, which we can call it the {\it residual cost}, and $\cJ^T_1(x;u(\cd))$ is roughly the long term average rate, which we hope to minimize.

\ms

Motivated by the above, we introduce the following so-called {\it probabilistic cell problem}.

\ms


\textbf{Problem (PC)}. For any $x\in\dbR^n$, determine a pair $(V(\cd),c)$
satisfying
\bel{eq:ergo_lq}\ba{ll}
\ns\ds c=\inf_{u(\cd)\in\sU}\limsup_{T\to\i}\frac{1}{T}
\int_0^T \mathbb E [f(X(t),u(t))]dt,\\
\ns\ds V(x)= \inf_{u(\cd) \in \sU}  \limsup_{T\to\i}\int_0^T\dbE[f(X(t),u(t))-c]dt.
\ea\ee
If there exists $\bar u(\cdot) \in \sU$ together with its corresponding state process $\bar X(\cdot)$
that attains the above two infimums,
the pair $(\bar X(\cdot), \bar u(\cdot))$ is referred to as the {\it optimal pair}.

\ms

In a standard control problem,
the objective is to find an optimal pair
$(\bar X(\cd),\bar u(\cd))$ to minimize the cost functional
to get the {\it value function} $V(\cd)$;
rather, Problem (PC) is somehow a two-objective problem:
seek a pair $(\bar X(\cdot),\bar u(\cdot))$ to
minimize the long term average rate and,
at the same time, to minimize the long term residual cost as well. We see that Problem (PC) admits a more general setting than that of Problem (LQ)$^\i_0$, mainly allowing $b,\si,q$ and $r$ to be nonzero.

\ms

It is worth noting that, in such a problem, the choice of $u(\cd)$ does not need to guarantee the existence of $\ds\lim_{T\to\i}J^T(x;u(\cd))$, in other words, $\sU$ might not be a subset of $\sU_{ad}[0,\i)$ (see \rf{U_ad[0,i)}). We now describe this new control space $\sU$. To this end, let $\cP_2(\dbR^k)$ be the Wasserstein space of probability measures $\m$ on $\dbR^k$ satisfying $\int_{\dbR^k}|x|^2d\m(x)<\i$, endowed with $2$-Wasserstein metric $\cW_2(\cd\,, \cd)$ defined by
$$\cW_2(\m_1,\m_2)=\inf_{\pi\in\Pi(\m_1,\m_2)}\(\int_{\dbR^k}\int_{\dbR^k}|x-y|^2
d\pi(x, y)\)^\frac{1}{2},$$
where $\Pi(\m_1,\m_2)$ is the collection of all probability measures on $\dbR^k\times\dbR^k$ with its marginals agreeing with $\m_1$ and $\m_2$, respectively. Denote
$$\int_{\dbR^k}\phi(x)\n(dx)=\lan\phi,\n\ran,$$
for all continuous function $\phi$ valued in $\dbR^k$ and all $\n\in\cP_2(\dbR^k)$. 

We define $\sU$ as the set of all processes $u(\cdot)$ such that for all $x\in\dbR^n$, the state process $X(\cd)\equiv X(\cd\,;x,u(\cd))$ of the solution to equation \rf{state}
satisfying
\bel{EX(t)|^2}
\dbE\big[|X(t)|^4+|u(t)|^4\big]<\i, \ee
and there exists a measure $\m_\i\in\cP_2(\dbR^n\times\dbR^m)$ such that
\bel{W_2}\lim_{t\to\i}\cW_2(\cL(X(t),u(t)),\m_\i)=0,\ee
where $\cL(X(t),u(t))$  is the law of $(X(t),u(t))$.

\ms

We will construct a $u(\cd)\in\sU$ below in Theorem \ref{t:verify}
to get the non-emptiness of this set. The following lemma will be useful in the proof.

\bl{l:convx} \sl Let $u(\cd)\in\sU$ such that $\cL(X(t),u(t))$
convergent to $\m_\i$ in 2-Wasserstein distance. Then for any continuous function $\phi:\dbR^n\times\dbR^m\to\dbR$ with a quadratic growth
$$|\phi(x,u)|\les K(1+|x|^2+|u|^2),\qq\forall(x,u)\in\dbR^n\times\dbR^m,$$\
it holds
$$\lim_{t\to\i}\dbE[\phi(X(t),u(t))]=\lan\phi,\m_\i\ran.$$

\el

\it Proof. \rm Since $\m_t\equiv\cL(X(t),u(t))$ converges to some $\m_\i\in\cP_2(\dbR^n\times\dbR^m)$ in 2-Wasserstein distance,
by Skorohod representation theorem (see \cite{Billingsley 1995}), one can find another possibly different probability space on which stochastic processes $(\wt X(\cd),\wt u(\cd))$ are defined with $\cL(\wt X(t),\wt u(t))=\m_t$ for all $t>0$, $\cL(\wt X_\i,\wt u_\i)=\m_\i$, and $(\wt X(t),\wt u(t))\to(\wt X_\i,\wt u_\i)$ as $t\to\i$ almost surely. Hence, by the fact of
$$\phi(\wt X(t),\wt u(t))\les K\big(1+|\wt X(t)|^2+|\wt u(t)|^2\big)$$
and
$$\lim_{t\to\i}\dbE\big[K\big(1+|\wt X(t)|^2+|\wt u(t)|^2\big)\big]=\dbE \big[K\big(1+|\wt X_\i|^2+|\wt u_\i|^2\big)\big],$$
one can apply the dominated convergence theorem to $(\wt X(t),\wt u(t))$ and obtain
$$\lim_{t\to\i}\dbE[\phi(X(t),u(t))]=\lim_{t\to\i}\dbE \big[\phi(\wt X(t),\wt u(t)) \big]
=\dbE \Big[\lim_{t\to\i}\phi(\wt X(t),\wt u(t)) \Big]=\lan\phi,\m_\i\ran.$$
This proves our conclusion.
\endpf

\ms

Under the conditions of Lemma \ref{l:convx}, for any $\phi(\cd)$ as above,  we have the following result:
\bel{5.9}\lim_{T\to\i}\frac{1}{T}\int_0^T\dbE[\phi(X(t),u(t))]dt=\lan\phi,\m_\i\ran.\ee

\ms

\subsection{Verification of Problem (PC) for LQ setting}
\label{s:verificantion for LQG setting}

In this subsection, we provide a complete characterization of
Problem (PC) in the LQ setting by the help of Theorem \ref{p:chara_ergodic_cost1} from Problem (C). We first introduce
the following SDE, whose solution $\bar X(\cd)$ turns out to be the optimal state of Problem (PC) (see a proof in Theorem \ref{t:verify})
\bel{bar X}\left\{\2n\ba{ll}
\ns\ds d\bar X(t)=\big[\bar A\bar X(t)+\bar b\big]dt
+\big[\bar C\bar X(t)+\bar\si\big]dW(t),\qq t\ges0,\\
\ns\ds\bar X(0)=x,\ea\right.\ee
where $\bar A,\bar C,\bar b,\bar\si$ are given in \rf{bar A bar C}. Thus, the above is the same as \rf{BX-SDE} with a possibly different initial condition.

\ms

Denote $\pi_\#$ as the pushforward measure of the projection $\pi(x,u)=x$, i.e.,
$$\pi_\#\m(B)=\m(\pi^{-1}(B)),\qq\forall B\in\cB(\dbR^n),\q\m\in\cP_2(\dbR^n\times\dbR^m).$$
The following theorem offers a comprehensive characterization of Problem (PC) within the LQ framework, facilitated by the results from Problem (C) outlined in Theorem \ref{p:chara_ergodic_cost1}.

\bt{t:verify} \sl Suppose {\rm(H1)--(H2)} hold. Let $(V(\cd),c_0)\in C^2(\dbR^n)\times\dbR$ be the solution to Problem {\rm(C)} of the form
$$V(x)=\frac{1}{2}\lan Px,x\ran+\lan p,x\ran,$$
where $P\in\dbS^n_{++}$ is the stabilizing solution of the algebraic Riccati equation \rf{ARE}, and $p$ is given by \rf{p}.
Let $\bar X(\cd)$ be the solution to the SDE \rf{bar X},
and $\bar u(\cd)$ be given in the form of
\bel{eq:ustar}\bar u(t)=\h u(\bar X(t),V_x(\bar X(t)),V_{xx}(\bar X(t))).\ee
%
where $\h u:\dbR^n\times\dbR^n\times\dbR^{n\times n}\mapsto\dbR^m$ is as defined in \eqref{hu}, and is rewritten below for clarity:
$$\h u(x,\Bp,\BP)=-(R\1n+\1n D^\top\BP D)^{-1}[(D^\top\BP C+S)x+D^\top\BP\si+B^\top\Bp+r].$$
Then, the following hold:

\ms

{\rm(i)} The distribution of the process $(\bar X(t),\bar u(t))$ converges to some $\bar\m_\i\in\cP_2(\dbR^n\times\dbR^m)$ in 2-Wasserstein distance, as $t\to \i$;

\ms

{\rm(ii)} Constant $c_0$ is uniquely determined by the following, called {\it ergodic constant}
\bel{eq:ergo_c0}c_0=\lim_{T\to\i}\frac{1}{T}\int_0^T\dbE\big[f(\bar X(t),\bar u
(t))\big]dt;\ee

\ms

{\rm(iii)} The 4-tuple $\{V(\cd)-\lan V,\pi_\#\bar\m_\i\ran,c_0,\bar X(\cd),\bar u(\cd)\}$ solves Problem {\rm (PC)}.

\ms

\et

\begin{remark}
    \rm Note that, without loss of generality, we choose $V(\cd)$ in the form of equation \eqref{v} given in Theorem \ref{p:chara_ergodic_cost1} with $p_0=0$. This is because for any constant $p_0$, the solution $\{V(\cd)-\lan V,\pi_\#\bar\m_\i\ran,c_0,\bar X(\cd),\bar u(\cd)\}$ to Problem (PC) does not depend on $p_0$.

\end{remark}

\begin{remark}
It is also worth noting that Theorem \ref{t:verify} (i) shows that the measure $\bar\mu_\infty$ is flow-invariant with respect to $(\bar X(\cdot), \bar u(\cdot))$. By Definition 7.11 in \cite{Tran 2021}, $\bar \mu_\infty$ can be called a {\it Mather measure}. The closure of the union of the supports of all such Mather measures is referred to as the {\it Mather set}. This observation establishes a connection between our Problem (PC) and the weak Kolmogorov--Arnold--Moser (KAM) theory, first developed by Fathi \cite{Fathi 2003} and Mather \cite{Mather 1991}. The weak KAM theory provides a link between deterministic/stochastic control problems and their corresponding cell problems, offering a representation of the optimal ergodic cost. Further details on weak KAM theory can be found in \cite{Evans 2008, Tran 2021}.

\end{remark}

Before proving the theorem, we present two lemmas. The first one is stated as follows.

\bl{l:moment_boundness} \sl Let {\rm(H1)--(H2)} hold. Then, there exists a constant $K>0$, independent of $t$, such that the solution $\bar X(\cd)$ to the SDE \rf{bar X} satisfies, for all $t\ges0$,
$$\dbE\big[|\bar X(t)|^2\big]\les K \q \mbox{and}\ \ \   \mathbb E\big[|\bar X(t)|^4\big]<\i.$$
\el

\it Proof. \rm First of all, by using (D.5) of \cite{FS06},
we have the following estimate (for any $p\ges 2$)
\bel{E|X|}\ba{ll}
\ns\ds\dbE\[\sup_{t\in[0,T]}|\bar X(t)|^p\]\les K |x|^p + KT^{\frac{p}{2}-1}e^{KT} \dbE\[\int_0^T (|x|^p + |\bar b|^p + |\bar\si|^p) dt\]\\
\ns\ds\qq\qq\qq\q~\les K\big(|x|^p+(|x|^p + |\bar b|^p + |\bar \sigma|^p) T^\frac{p}{2}e^{KT}\big).\ea\ee
Thus, the second result is obtained by choosing $p = 4$. Next, since $[\bar A,\bar C]$ is $L^2$-exponentially stable, by Lemma \ref{l:exponential_stable}, the following Lyapunov equation has a unique solution $P\in\dbS^n_{++}$:
$$P\bar A+\bar A^\top P+\bar C^\top P \bar C+I=0.$$
Let $\beta > 0$ be the largest eigenvalue of $P$. Now, applying It\^o's formula to $t\mapsto\lan P\bar X(t),\bar X(t)\ran$, we have
\begin{equation*}
\begin{aligned}
& \frac{d}{dt} \dbE[\lan P\bar X(t),\bar X(t)\ran] \\
= \ & \dbE[\lan(P\bar A+\bar A^\top P+\bar C^\top P\bar C)\bar X(t),\bar X(t)\ran + 2 \lan \bar X(t), P \bar b+ \bar C^\top P \bar \sigma \ran + \lan P \bar \sigma, \bar \sigma \ran] \\
\les \ & 
 \dbE[- |\bar X(t)|^2 + 2 |P \bar b+ \bar C^\top P \bar \sigma| |\bar X(t)|+\lan P \bar \sigma, \bar \sigma \ran] \\
\les \ & -\frac{1}{2} \dbE[ |\bar X(t)|^2] + 2 |P \bar b+ \bar C^\top P \bar \sigma|^2 + \lan P \bar \sigma, \bar \sigma \ran \\
\les \ & -\frac{1}{2 \beta} \dbE[\lan P\bar X(t),\bar X(t)\ran] + 2 |P \bar b+ \bar C^\top P \bar \sigma|^2 + \lan P \bar \sigma, \bar \sigma \ran.
\end{aligned}
\end{equation*}
Then, denoting $\delta = \frac{1}{2 \beta}$ and applying Gr\"onwall's inequality, we get the following estimation
$$\dbE[\lan P\bar X(t),\bar X(t)\ran] \les e^{-\d t}\lan Px,x\ran + \frac{2 |P \bar b+ \bar C^\top P \bar \sigma|^2 + \lan P \bar \sigma, \bar \sigma \ran}{\d} (1 - e^{-\d t})$$
for all $t \ges 0$.
Hence, by the positivity of $P$, we obtain our conclusion.
\endpf

\ms

The second lemma states the asymptotic behavior of the distribution of some linear nonhomogeneous uncontrolled SDE.
\ms

\bl{l:verify} \sl If system $[A,C]$ is $L^2$-exponentially stable, then the solution $X(\cd)$ to the following SDE
\bel{eq:uncon}dX(t)=[AX(t)+b]dt+[CX(t)+\si]dW(t),\qq X(0)=x,\ee
allowing $b,\si\ne0$, has an invariant measure $\n_\i$ in the $2$-Wasserstein space $\cP_2(\dbR^n)$ such that
$$\cW_2(\cL(X(t)),\n_\i)\to0,\qq\hb{as }t\to\i.$$
%


\el

\it Proof. \rm Let $\n_t=\cL(X(t))$. We claim that it is Cauchy in $\cP_2(\dbR^n)$. First, let $\F(\cd)$ be the fundamental matrix of $[A,C]$, i.e.,
$$d\F(t)=A\F(t)dt+C\F(t)dW(t),\qq\F(0)=I.$$
It is known that $\F(t)$ is invertible for all $t\ges0$. Then, the explicit solution of the nonhomogeneous SDE \rf{eq:uncon} in terms of the above fundamental solution is given by
$$X(t)=\F(t)x+\F(t)\int_0^t\F(r)^{-1}(b-C\si)dr+\F(t)\int_0^t\F(r)^{-1}\si d W(r).$$
Denote that $\F(s,t)=\F(t)\F(s)^{-1}$. The nonhomogeneous SDE \rf{eq:uncon} satisfies the following representation exhibiting the Markov property:
\bel{eq:expli_soltion_st}\ba{ll}
\ns\ds X(t+s)=\F(s,s+t)X(s)+\F(s,s+t)\int_s^{s+t}\F(s,r)^{-1}(b-C\si)dr\\
\ns\ds\qq\qq\qq+\F(s,s+t)\int_s^{s+t}\F(s,r)^{-1}\si dW(r).\ea\ee
Next, let $\h\F(\cd)$ be the solution to the following
$$d\h\F(t)=A\h\F(t)dt+C\h\F(t)d\h W(t),\qq \h\F(0)=I,$$
where $\h W(\cd)$ is another standard Brownian motion, independent of $W(\cd)$. It is also known that $\h\F(t)^{-1}$ exists. We now define another process $\h X^s(\cd)$ as the solution to the following nonhomogeneous equation
$$\h X^s(t)=\h\F(t)X(s)+\h\F(t)\int_0^t\h\F(r)^{-1}(b-C\si)dr+\h\F(t)\int_0^t \h\F(r)^{-1}\si d\h W(r).$$
Then, $\h X^s(t)$ and $X(t+s)$ have the same law. Similarly, we have a distribution copy of $X(t)$ by $\h X^0(t)$ given by
$$\h X^0(t)=\h\F(t)x+\h\F(t)\int_0^t\h\F(r)^{-1}(b-C\si)dr+\h\F(t)\int_0^t \h\F(r)^{-1}\si d\h W(r).$$
Therefore, by the definition of $L^2$-exponentially stability in Lemma \ref{l:exponential_stable} (i) and the similar argument about the $t$-uniform boundness of second moment in Lemma \ref{l:moment_boundness} , we have
\begin{equation*}
\begin{aligned}
\cW_2^2(\nu_{t+s},\nu_t) &= \cW_2^2 \big(\cL\big(\h X^s(t)\big),\cL\big(\h X^0(t)\big)\big)\les\dbE\big[\big|\h X^s(t)-\h X^0(s)\big|^2\big]\\
&=\dbE\big[\big|\h\Phi(t)(X(s)-x)\big|^2\big]\les Ke^{-\l t}\dbE \big[|X(s)-x|^2 \big]\les Ke^{-\l t},
\end{aligned}
\end{equation*}
where $K,\l>0$ are independent of $s$ and $t$. This implies that $\nu_t$ is a Cauchy sequence in $\cP_2(\dbR^n)$.

\endpf

\ms

Now, we are ready to prove Theorem \ref{t:verify}.

\ms

\it Proof of Theorem \ref{t:verify}. \rm The first result is a straightforward result of Lemma \ref{l:verify}. For simplicity of the notations, we write
$$\h b(x,u)\equiv Ax+Bu+b,\qq\h\si(x,u)\equiv Cx+Du+\si.$$
Since $(V(\cd),c_0)\in C^2(\dbR^n)\times\dbR$ solves Problem (C), for a control $u(\cd)\in\sU$ and its associated state process $X(\cd) = X(\cd\ ; x, u(\cd))$, by It\^o's formula, we have
$$\ba{ll}
\ns\ds V(X(t))=V(x)+\int_0^t\big(\lan\h b(X(s),u(s)),V_x(X(s))\ran+\frac{1}{2}\lan V_{xx}(X(s))\h\si(X(s),u(s)),\h\si(X(s),u(s))\ran\big)ds\\
\ns\ds\qq\qq\qq+\int_0^t\lan\h\si(X(s),u(s)),V_x(X(s))\ran dW(s).\ea$$
Note that the quadratic variation of the last term is
\begin{equation*}
\begin{aligned}
\dbE\[\int_0^t|\lan\h\si(X(s),u(s)),V_x(X(s))\ran|^2ds\] &\les K\dbE\[ \int_0^t\big(1+|X(s)|^4+|u(s)|^4\big)ds\]\\
&\les K\(t+\dbE\[\int_0^t\big(|X(s)|^4+|u(s)|^4\big)ds\]\),
\end{aligned}
\end{equation*}
and it is finite by the condition $u(\cd)\in\sU$. So, fixing $t>0$, and taking expectation on both sides, we obtain
\bel{EV=V+}\ba{ll}
\ns\ds\dbE\big[V(X(t))\big]=V(x)+\dbE\[\int_0^t\big(\lan\h b(X(s),u(s)),V_x(X(s))\ran\\
\ns\ds\qq\qq\qq\qq+\frac{1}{2}\lan V_{xx}(X(s))\h\si(X(s),u(s)),\h\si(X(s),u(s)) \ran\big)ds\].\ea\ee
Problem (C) implies that, for all $(X(\cd),u(\cd))$ and $s\in[0,t]$,
\bel{dbH>c_0}\ba{ll}
\ns\ds\lan\h b(X(s),u(s)),V_x(X(s))\ran+\frac{1}{2}\lan V_{xx}(X(s))\h\si(X(s), u(s)),\h\si(X(s),u(s))\ran+f(X(s),u(s))\\
\ns\ds=\dbH(X(s),V_x(X(s)),V_{xx}(X(s)),u(s))\ges H(X(s),V_x(X(s)),V_{xx}(X(s)))=c_0.\ea\ee
Thus, combining \rf{EV=V+} and \rf{dbH>c_0}, we have
\begin{equation}
\label{eq:vx1}
\begin{aligned}
V(x) &= \dbE\big[V(X(t))\big]-\dbE\[\int_0^t\big(\lan\h b(X(s),u(s)),V_x(X(s))\ran\\
&\qq\qq+\frac{1}{2}\lan V_{xx}(X(s))\h\si(X(s),u(s)),\h\si(X(s),u(s)) \ran\big)ds\]\\
& \les\dbE\big[V(X(t))\big]+\dbE\[\int_0^t\big(f(X(s),u(s))-c_0\big)ds\].
\end{aligned}
\end{equation}
The inequality \rf{eq:vx1} holds for all $u(\cd)\in\sU$. Clearly, Lemmas \ref{l:moment_boundness} and \ref{l:verify}, together with equation \rf{eq:ustar} imply that $\bar u(\cd)\in\sU$. This, by the way, shows that $\sU\ne\varnothing$. Now, in Theorem \ref{p:chara_ergodic_cost1} by taking $p_0=0$, we know that the equality in \rf{eq:vx1} holds if $u(\cd)=\bar u(\cd)$. Since the value function $V(\cd)$ satisfies a quadratic growth condition, i.e., $|V(x)|\les K(1+|x|^2)$ for all $x\in\dbR^n$, by Lemma \ref{l:convx}, we have
\bel{EV}\ba{ll}
\ns\ds\lim_{t\to\i}\dbE\big[V(X(t))\big]=\lan V,\pi_\#\m_\i^u\ran,\qq\lim_{t\to\i}\dbE\big[f(X(t),u(t))\big]=\lan f,\m_\i^u\ran,\ea\ee
where $\m_\infty^u$ is the distribution limit of $(X(t),u(t))$
as $t\to\i$, for a given $u(\cd)\in\sU$.
Note that the superscript $u$ is used to
emphasize the dependence
on the control $u(\cd)$.

\ms

Dividing $t$ on both sides of \rf{eq:vx1} and then sending $t\to\i$, we get
\bel{eq:c0p}c_0\les\lim_{t\to\i}\frac{1}{t}\int_0^t\dbE\big[f(X(s),u(s))\big]ds=\lan f,\m_\i^u\ran.\ee
The above inequality holds for all $u(\cd)\in\sU$ and the equality holds if $u(\cd)=\bar u(\cd)$, with $\m_\i^{\bar u}=\bar\m_\i$. Thus, we conclude \rf{eq:ergo_c0}, or (ii) of  Theorem \ref{t:verify}.

\ms

Next, sending $t\to\i$ on both sides of \rf{eq:vx1}, we have
\bel{eq:vp1}V(x)\les\lan V,\pi_\#\m_\i^u\ran+\lim_{t\to\i}
\int_0^t\dbE\big[f(X(s),u(s))-c_0\big]ds.\ee
It leads to the following inequality: for any $u(\cd)\in\sU$,
\bel{eq:vp2}V(x)\les\lan V,\pi_\#\bar\m_\i\ran+\lim_{t\to\i}
\int_0^t\dbE\big[f(X(s),u(s))-c_0\big]ds.\ee
This is because the following: On one hand, due to Lemma \ref{l:convx} and \rf{eq:c0p},
$$\lan f,\m_\i^u\ran=\lim_{t\to\i}\dbE\big[f(X(t),u(t))\big]
=\lim_{T\to\i}\frac{1}{T}\int_0^T\dbE\big[f(X(t),u(t))\big]dt
\ges \lan f,\bar\m_\i\ran=c_0.$$
Therefore, it is enough to consider two cases:
If $\lan f,\m_\i^u\ran=\lan f, \bar \m_\i\ran$,
then \rf{eq:vp1} holds with $\bar\m_\i$ in place of $\m_\i^u$, which yields \eqref{eq:vp2};
If $\lan f,\m_\i^u\ran>\lan f,\bar\m_\i\ran=c_0$, then it implies that
$$\lim_{t \to \infty}\int_0^{t} \dbE\big[f(X(s),u(s))-c_0\big]ds=\infty,$$
which yields \eqref{eq:vp2} with the right-hand side being $\i$.

\ms

Moreover, the equality in \rf{eq:vp2} holds when $u(\cd)=\bar u(\cd)$.
Therefore, we rewrite it as
$$V(x)-\lan V,\pi_\#\bar\m_\i\ran=\inf_{u(\cd)\in\sU}\lim_{t\to\i}
\int_0^t\dbE\big[ f(X(s),u(s))-c_0\big]ds,$$
and (iii) of the theorem follows. \endpf


\section{Turnpike Property from Probabilistic Cell Problem}
\label{s:convergence}

As we indicated earlier, Problem (LQ)$^T$
enjoys the turnpike property (as $T\to\i$).
In obtaining such a property, the horizon independent process
$(\BX^*(\cd),\Bu^*(\cd))$ given explicitly by \rf{eq:tp_yong}
is used. See \cite{Sun-Yong 2024} for details.
On the other hand, Problem (PC) is clearly closely related to Problem (LQ)$^T$, for large $T>0$.
Hence, we expect that one might be able to use the solution of Problem (PC) to exhibit turnpike properties for Problem (LQ)$^T$. In this section, we are going to reveal such a possibility, and present the turnpike property for Problem (LQ)$^T$ using an alternative way.

\ms

Now we present our first main result of this section, a turnpike property between the optimal pair of Problem (LQ)$^T$ and those of Problem (PC).

\bt{t:convergence_of_optimal_path} \sl
Let {\rm(H1)--(H2)} hold. Let $(X^T(\cd), u^T(\cd))$ be the optimal pair of Problem (LQ)$^T$ with initial state $x \in \mathbb R^n$ and $(\bar X(\cd), \bar u(\cd))$ be the optimal pair of Problem (PC) with initial state $\bar x \in \mathbb R^n$. Then, there exist some constants $\lambda, K > 0$, independent of $T$, such that
\begin{equation}
\label{eq:strong_turnpike}
\mathbb E \big[|X^{T}(t) - \bar X (t)|^2
+ |u^{T} (t) - \bar u (t)|^2
\big] \les K \big(e^{-\lambda t} + e^{-\lambda(T - t)} \big), \qq \forall t \in [0, T].
\end{equation}
\et

\it Proof. \rm Let $(\bar X(\cd),\bar u(\cd))$
be optimal pair of Problem (PC) and
$$\h X(t)=X^T(t)-\bar X(t),\qq t\in[0,T].$$
Then the dynamic of $\h X(\cdot)$ is given by the following SDE
\begin{equation*}
\begin{cases}
d\h X(t)=\big((A+B\Th^T(t))\h X(t)+B(\Th^T(t)-\bar\Th)\bar X(t)+B(\th^T(t)-\bar\th)\big)dt\\
\qq\qq\q+\big((C+D\Th^T(t))\h X(t)+D(\Th^T(t)-\bar\Th)\bar X(t)+D(\th^T(t)-\bar\th)\big)dW(t),\\
\h X(0)= x - \bar x.
\end{cases}
\end{equation*}
Let $P$ be the positive definite solution to the algebraic Riccati equation \eqref{ARE}. The It\^o's formula implies that
$$\ba{ll}
\ns\ds\frac{d}{dt}\dbE\big[\lan P\h X(t),\h X(t)\ran\big]\\
\ns\ds=\dbE\big[2\lan P\h X(t),(A+B\Th^T(t))\h X(t)+B(\Th^T(t)-\bar\Th) \bar X(t)+B(\th^T(t)-\bar\th)\ran\\
\ns\ds\qq+\lan P[(C+D\Th^T(t))\h X(t)+D(\Th^T(t)-\bar\Th)\bar X(t)+D(\th^T(t)-\bar\th)],\\
\ns\ds\qq\qq[(C+D\Th^T(t))\h X(t)+D(\Th^T(t)-\bar\Th)\bar X(t)+D(\th^T(t)-\bar \th)]\ran\big]\\
\ns\ds=\dbE\big[\lan[P(A+B\Th^T(t))+(A+B\Th^T(t))^\top P+(C+D\Th^T(t))^\top P (C+D\Th^T(t))]\h X(t),\h X(t)\ran\\
\ns\ds\qq+2\lan\h X(t),\h H(t)\ran\\
\ns\ds\qq+\lan P[D(\Th^T(t)-\bar\Th)\bar X(t)+D(\th^T(t)-\bar\th)],[D(\Th^T(t)-\bar\Th)\bar X(t)+D(\th^T(t)-\bar \th)]\ran\big],\ea$$
where, for all $t \in [0, T]$,
$$\ba{ll}
\ns\ds\h H(t)\equiv P[B(\Th^T(t)-\bar\Th)\bar X(t)+B(\th^T(t)-\bar\th)]\\
\ns\ds\qq\qq+(C+D\Th^T(t))^\top P[D(\Th^T(t)-\bar\Th)\bar X(t)+D(\th^T(t)-\bar\th)].\ea$$
By Theorem \ref{approximation} and Lemma \ref{l:moment_boundness}, there exist constants $K,\l>0$, independent of $T$, such that for all $t \in [0, T]$,
$$\ba{ll}
\ns\ds\dbE\big[\lan P[D(\Th^T(t)-\bar\Th)\bar X(t)+D(\th^T(t)-\bar\th)], [D(\Th^T(t)-\bar\Th)\bar X(t)+D(\th^T(t)-\bar\th)]\ran\big]\\
\ns\ds\les K\dbE\big[|D(\Th^T(t)-\bar\Th)\bar X(t)|^2+|D(\th^T(t)-\bar \th)|^2\big]\\
\ns\ds\les K\big(|\Th^T(t)-\bar\Th|^2\dbE\big[|\bar X(t)|^2\big]+|\th^T(t)-\bar \th|^2\big) \les Ke^{-\l(T-t)}.\ea$$
Next, we observe the facts that
$$\ba{ll}
\ns\ds P(A+B\Th^T(t))+(A+B\Th^T(t))^\top P
+(C+D\Th^T(t))^\top P (C+D\Th^T(t))\\
\ns\ds=P(A+B\bar\Th)+(A+B\bar\Th)^\top P+(C+D\bar\Th)^\top P(C+D\bar\Th)+ PB(\Th^T(t)-\bar\Th)\\
\ns\ds\qq+(B(\Th^T(t)-\bar\Th))^\top P+(D(\Th^T(t)-\bar\Th))^\top P (C+D\bar\Th)+(C+D\Th^T(t))^\top PD(\Th^T(t)-\bar\Th),\ea$$
and
$$P(A+B\bar\Th)+(A+B\bar\Th)^\top P+
(C+D\bar\Th)^\top P(C+D\bar\Th)=
-(Q + S^\top\bar\Th+\bar\Th^\top S+\bar\Th^\top R\bar\Th)<0.$$
Theorem \ref{approximation} shows that, for all $t \in [0, T]$,
\begin{equation*}
\begin{aligned}
& |PB (\Th^T(t) - \bar\Th) + (B (\Th^T(t) - \bar\Th))^\top P + (D (\Th^T(t) - \bar\Th))^\top P (C+D\bar\Th) \\
& \hspace{0.5in} + (C+D\Th^T(t))^\top P D (\Th^T(t) - \bar\Th) | \les Ke^{-\lambda(T-t)}.
\end{aligned}
\end{equation*}
Let $\beta_1 > 0$ and $\beta_2 > 0$ be the smallest eigenvalues of $P$ and $Q + S^\top \bar \Th + \bar \Th^\top S + \bar \Th^\top R \bar \Th$ respectively, and $\h \beta > 0$ is the largest eigenvalue of $P$. Then, it is clear that
\begin{equation}
\label{temp}
\beta_1 \dbE\big[|\h X(t)|^2\big] \les \dbE \big[\lan P \h X(t), \h X(t)\ran \big] \les \h \beta \dbE\big[|\h X(t)|^2\big]
\end{equation}
for all $t \in [0, T]$.
The above results and Young's inequality imply that
\begin{equation*}
\begin{aligned}
\frac{d}{dt} \dbE \big[\lan P \h X(t), \h X(t)\ran \big] & \les \dbE \big[\big(- \beta_2 + Ke^{-\lambda(T-t)} \big) |\h X(t)|^2 + 2 |\h H(t)| |\h X(t)| + Ke^{-\lambda(T-t)} \big] \\
& \les \dbE \Big[ \(- \frac{\beta_2}{2} + Ke^{-\lambda(T-t)} \) |\h X(t)|^2 + \frac{2}{\beta_2} |\h H(t)|^2 + Ke^{-\lambda(T-t)} \Big]
\end{aligned}
\end{equation*}
for all $t \in [0, T]$. Note that, by Theorem \ref{approximation} again, there exist some potentially distinct constants $K, \lambda > 0$, independent of $T$, such that
$$\dbE \big[ |\h H(t)|^2 \big] \les Ke^{-\lambda(T-t)}, \qq \forall t \in [0, T].$$
Hence, by the inequality \eqref{temp}, there exists a constant $K_1 > 0$ such that
$$\frac{d}{dt} \dbE \big[\lan P \h X(t), \h X(t)\ran \big] \les K_1 \(- \frac{\beta_2}{2} + K e^{-\lambda(T-t)} \) \dbE \big[\lan P \h X(t), \h X(t)\ran\big] + K e^{-\lambda(T-t)}$$
for all $t \in [0, T]$. Denote that $g(t) := K_1(- \frac{\beta_2}{2} + Ke^{-\lambda(T-t)})$ for all $t \in [0, T]$. Then, it is clear that, for all $0 \les s \les t < T$, the following estimate holds:
$$\exp \Big\{\int_s^t g(r) dr\Big\} = \exp \Big\{\frac{K_1K}{\lambda} e^{-\lambda T}(e^{\lambda t} - e^{\lambda s}) - \frac{K_1 \beta_2}{2}(t-s) \Big\} \les \exp \Big\{\frac{K_1K}{\lambda} - \frac{K_1 \beta_2}{2}(t-s) \Big\}.$$
Thus, we have
$$\int_0^t e^{\int_s^t g(r) dr} K e^{-\lambda (T-s)} ds \les K e^{-\lambda (T-t)}, \quad \forall t \in [0, T],$$
with possible different constants $K, \lambda > 0$.
It follows that
\begin{equation}
\label{temp2}
\begin{aligned}
\dbE \big[\lan P \h X(t), \h X(t)\ran \big] & \les \lan P(x - \bar x), x - \bar x \ran e^{\int_0^t g(r) dr} + \int_0^t e^{\int_s^t g(r) dr} K e^{-\lambda (T-s)} ds \\
& \les K \big(e^{-\lambda t} + e^{-\lambda(T-t)} \big)
\end{aligned}
\end{equation}
for all $t \in [0, T]$ with potentially distinct constants $K > 0$ and $\lambda > 0$. Since $\beta_1$ is the smallest eigenvalues of $P$, then we obtain the desired result that
$$\dbE\big[|\h X(t)|^2\big] \les \frac{1}{\beta_1} \dbE \big[\lan P \h X(t), \h X(t)\ran \big] \les  K \big(e^{-\lambda t} + e^{-\lambda(T - t)} \big), \qq \forall t \in [0, T]$$
for some constants $K > 0$ and $\lambda > 0$ independent of $t$ and $T$.

From the above estimate and the result in Lemma \ref{l:moment_boundness}, it is clear that $\mathbb E[|X^T(t)|^2] \les K$ for all $t \in [0, T]$ for some $K > 0$. Recall that
\begin{equation*}
u^T(t) = \Th^T(t) X^T(t) + \theta^T(t),
\qq \bar u(t) = \bar \Theta \bar X (t) + \bar \theta,
\end{equation*}
then we deduce the following estimate
\begin{equation*}
\begin{aligned}
\dbE \big[ |u^T(t) - \bar u (t)|^2 \big]
& \les K \dbE \big[ |(\Th^T(t) - \bar \Th) \bar X^T(t) |^2 +
|\bar \Th(X^T(t) - \bar X (t)) |^2 + |\theta^T(t) - \bar \theta|^2 \big] \\
& \les K \big(e^{-\lambda t} + e^{-\lambda(T - t)} \big)
\end{aligned}
\end{equation*}
for all $t \in [0, T]$ by using the results in Theorem \ref{approximation}. This completes the proof.
\endpf

\ms

\begin{remark} \rm
\label{initial}
In this result, we assume different initial conditions, $x$ and $\bar{x}$, for Problem (LQ)$^T$ and Problem (PC), respectively. An interesting observation arises from equation \eqref{temp2}. When $x = \bar{x}$, the stochastic turnpike property only returns a one-sided result as
$$\mathbb E \big[|X^{T}(t) - \bar X (t) |^2
+ |u^{T} (t) - \bar u (t) |^2
\big] \les K e^{-\lambda(T - t)}, \qq \forall t \in [0, T].$$
We can intuitively understand this as: when the dynamics of these two problems start from the same initial value, the optimal paths of the state processes and the controls coincide at the beginning and will deviate only as they approach the terminal time.


\end{remark}

\ms

This main result provides the following straightforward corollary for the mean-square turnpike property.

\begin{corollary} \sl Let {\rm(H1)--(H2)} hold. Then,
$$\lim_{T\to\i}\frac{1}{T}\int_0^T\dbE\big[|X^T(t)-\bar X(t)|^2+|u^T(t)-\bar u(t)|^2\big]dt=0.$$

\end{corollary}

The next corollary demonstrates that, for any initial
state $x$, the value function $V^{T}(x)$ of Problem (LQ)$^T$
converges, in the time-average sense,
to the optimal value $L(x^*, u^*)$
from the static optimization problem \eqref{S-static}. Its proof can be found in Corollary 3.4 
in \cite{Sun-Yong 2024}.

\begin{corollary} \sl
\label{c:consistence_static}
Let {\rm(H1)--(H2)} hold. Then, for all $x \in \mathbb R^n$,
\begin{equation}
\label{eq:convergence_average_value}
\lim_{T \to \infty} \frac{1}{T} V^{T} (x) = L(x^*, u^*).
\end{equation}
\end{corollary}

\ms

From the above result \eqref{eq:convergence_average_value} and the definition of ergodic cost in \eqref{eq:ecost1}, it is straightforward to see that the optimal value $L(x^*, u^*)$ from the static optimization problem \eqref{S-static}
is consistent with the ergodic cost $\bar c$.

Distinct from the aforementioned turnpike property elucidated by  \eqref{stoch-turnpike}, we also establish the turnpike behavior concerning the cost function:
$$\lim_{T\to\i}\frac{1}{T}J^T(x;\bar u(\cd))=\bar c\equiv\lim_{T\to\i}\frac{1}{T}J^T(x;u^T(\cd)).$$
In other words, achieving a near optimality in terms of the average cost over an extended period does not require calculating the optimal control $u^T(\cd)$ for every terminal time $T$.
Instead, one only needs to compute the optimal
control $\bar u(\cd)$ for the probabilistic cell problem (PC).

\ms

To accomplish this objective, it becomes imperative to
establish a relation between the cost function
$J^T(x;\bar u(\cd))$
with the $\bar u(\cd)$ from the probabilistic cell problem
(PC),
and the corresponding function $V^T(x) = J^T(x, u^T(\cd))$
originating from
the finite time stochastic control problem (LQ)$^T$.

\ms

Let us revisit the definition of the value function $V(x)$
for the probabilistic cell problem (PC), as given in \eqref{eq:ergo_lq}:
\begin{equation*}
V(x) \equiv \inf_{u(\cdot) \in \sU}  \lim_{T \to \infty}
\int_0^T  \mathbb{E} \big[f\left(X(t), u(t) \right) - c_0 \big] d t.
\end{equation*}
Here, the optimal control and optimal path of the
probabilistic cell problem are denoted by
$\bar u(\cdot)$ and $\bar X(\cdot)$ respectively, which are characterized by Theorem \ref{t:verify}.
It's worth noting that the constant $c$ in the above definition is replaced by $c_0$ from the cell problem corresponding to the positive definite $P$ according to Theorem \ref{t:verify}.

The second main result of this section is the turnpike property between the cost functional on the finite-horizon evaluated at the optimal control obtained from Problem (LQ)$^T$ and Problem (PC).

\begin{theorem} \sl
\label{t:convergence_of_value_function}
Let {\rm(H1)--(H2)} hold. Let $J^{T}(x;\bar u(\cd))$
be the cost functional evaluated along the optimal control
$\bar u(\cd)$ of Problem (PC)  on the finite horizon $[0, T]$, i.e.,
$$J^T(x;\bar u(\cd))\equiv\dbE\Big[\int_0^Tf(\bar X(s),\bar u(s))ds\Big].$$
For all $x \in \mathbb R^n$, the following estimation holds:
$$0\les J^T(x;\bar u(\cd))-V^T(x)=O(1).$$
Moreover, the constant $c_0$ of Problem (C) matches with the ergodic cost $\bar c$ defined via \eqref{eq:ecost1},
i.e., for all $x \in \mathbb R^n$,
$$c_0=\bar c=\lim_{T\to \infty} \frac{1}{T} V^{T}(x).$$
\end{theorem}

\begin{remark} \rm
\label{TurnpikeJ}
From the results of Theorem \ref{t:convergence_of_value_function}, we could establish the
turnpike behavior in terms of the average cost function straightforwardly. Note that for all $x \in \mathbb R^n$,
$V^{T}(x) \les J^{T}(x;\bar u(\cd)) = V^{T}(x) + O(1)$, it follows that
\begin{equation}
\label{eq:tpk3}
\lim_{T\to \infty} \frac{1}{T} J^T(x; \bar u(\cd)) = \lim_{T\to \infty} \frac{1}{T} V^{T}(x) = \bar c.
\end{equation}
\end{remark}

Finally, we prove Theorem \ref{t:convergence_of_value_function}.

\ms

\it Proof of Theorem \ref{t:convergence_of_value_function}. \rm
By definition, for all $x \in \mathbb R^n$, it is clear that
$$\big|J^T(x; \bar u(\cd)) -
J^T(x; u^T(\cd))\big| \les \dbE\[\int_0^T \big|f(\bar X(t),\bar u(t))-f(X^T(t),u^T(t)) \big| dt\].$$
Then, using H\"older's inequality, we have
$$\ba{ll}
\ns\ds\dbE\Big[\int_0^T \big|f(\bar X(t),\bar u(t))-f(X^T(t),u^T(t)) \big|dt \Big]\\
\ns\ds=\dbE\Big[\int_0^T\Big|\frac{1}{2}\lan Q\bar X(t),\bar X(t)\ran+\lan S \bar X(t),\bar u(t)\ran+\frac{1}{2}\lan R\bar u(t),\bar u(t)\ran+\lan q, \bar X(t)\ran+\lan r,\bar u(t)\ran\\
\ns\ds\qq-\frac{1}{2}\lan Q  X^T(t),X^T(t)\ran\1n-\1n\lan S X^T(t),u^T(t)\ran\1n-\1n\frac{1}{2}\lan Ru^T(t),u^T(t)\ran\1n-\1n\lan q,X^T(t)\ran\1n-\1n\lan r,u^T(t)\ran\Big|dt\Big] \\
\ns\ds=\dbE\Big[\int_0^T\Big|\frac{1}{2}\lan Q(\bar X(t)-X^T(t)),\bar X(t)\ran -\frac{1}{2}\lan Q(X^T(t)-\bar X(t)),X^T(t)\ran\\
\ns\ds\qq+\lan S(\bar X(t)-X^T(t)),\bar u(t)\ran-\lan SX^T(t),u^T(t)-\bar u(t)\ran\\
\ns\ds\qq+\frac{1}{2}\lan R(\bar u(t)-u^T(t)),\bar u(t)\ran-\frac{1}{2}\lan R (u^T(t)-\bar u(t)),u^T(t)\ran\\
\ns\ds\qq+\lan q,\bar X(t)-X^T(t)\ran+\lan r,\bar u(t)-u^T(t)\ran \Big|dt \Big] \\
\ns\ds\les K\int_0^T\3n\Big\{\Big[\big(\dbE[|\bar X(t)|^2]\big)^\frac{1}{2}+\big(\dbE[|X^T(t)|^2]\big)^\frac{1}{2}+\big(\dbE[|\bar u(t)|^2]\big)^\frac{1}{2}+|q|\Big]\big(\dbE[|\bar X(t)-X^T(t)|^2]\big)^\frac{1}{2}\\
\ns\ds\qq+\Big[\big(\dbE[|X^T(t)|^2]\big)^\frac{1}{2}+\big(\dbE[|u^T(t)|^2]\big)^\frac{1}{2}+ \big(\dbE[|\bar u(t)|^2]\big)^\frac{1}{2}+|r|\Big]\big(\dbE[|\bar u(t)-u^T(t)|^2]\big)^\frac{1}{2} \Big\}dt \\
\ns\ds\les K \int_0^{T} \big(e^{-\lambda t} + e^{-\lambda(T - t)} \big) dt \les K \ea$$
from the results of Lemma \ref{l:moment_boundness}, Theorem \ref{t:convergence_of_optimal_path}, and the fact that
$$\mathbb E \big[|X^{T}(t)|^2 + | u^{T}(t)|^2 \big] \les K, \qq \forall t \in [0, T].$$
Therefore, we obtain the desired estimation that
$$\big|J^T(x;\bar u(\cd))-J^T(x;u^T(\cd)) \big|\les K$$
for all $x \in \mathbb R^n$. Moreover, by the definition of the value function $V^{T}(x)\equiv J^{T} \left(x; u^T(\cdot) \right)$ from Problem (LQ)$^T$, it is straightforward that for all $x \in \mathbb R^n$,
$0 \les J^T(x; \bar u(\cdot)) - J^T(x; u^T(\cdot))$.
Therefore, for all $x \in \mathbb R^n$, the following estimation holds:
$$0\les J^T(x;\bar u(\cd))-V^T(x)=O(1).$$
\endpf


\end{document}